\documentclass[graybox]{svmult}

\usepackage{type1cm}        
\usepackage{makeidx}         
\usepackage{graphicx}        
\usepackage{multicol}        
\usepackage[bottom]{footmisc}

\usepackage{newtxtext}       %
\usepackage{newtxmath}       

\usepackage{bm}
\usepackage{layouts}

\makeindex             
                       
\begin{document}

\title*{THU-Splines: Highly Localized Refinement on Smooth Unstructured Splines}
\titlerunning{THU-Splines}
\author{Xiaodong Wei}
\institute{Xiaodong Wei \at École Polytechnique Fédérale de Lausanne (EPFL), 1015 Lausanne, Switzerland, \email{xiaodong.wei@epfl.ch}}
%
%
\maketitle

\abstract{We present a novel method named \emph{truncated hierarchical unstructured splines} (THU-splines) that supports both local $h$-refinement and unstructured quadrilateral meshes. In a THU-spline construction, an unstructured quadrilateral mesh is taken as the input control mesh, where the degenerated-patch method \cite{ref:reif97} is adopted in irregular regions to define $C^1$-continuous bicubic splines, whereas regular regions only involve $C^2$ B-splines. Irregular regions are then smoothly joined with regular regions through the truncation mechanism \cite{ref:wei18}, leading to a globally smooth spline construction. Subsequently, local refinement is performed following the truncated hierarchical B-spline construction \cite{ref:giannelli12} to achieve a flexible refinement without propagating to unanticipated regions. Challenges lie in refining transition regions where a mixed types of splines play a role. THU-spline basis functions are globally $C^1$-continuous and are non-negative everywhere except near extraordinary vertices, where slight negativity is inevitable to retain refinability of the spline functions defined using the degenerated-patch method. Such functions also have a finite representation that can be easily integrated with existing finite element or isogeometric codes through B\'{e}zier extraction.}

\section{Introduction}
\label{sec:intro}

Isogeometric analysis (IGA) was proposed to tightly integrate computer-aided design (CAD) with engineering simulation via employing the same spline-based basis in both sides \cite{ref:hughes05, ref:cottrell09}. Local refinement and extraordinary vertices\footnote{An extraordinary vertex is an interior vertex shared by other than four quadrilateral faces.} (EVs) have been two of the most intensively studied topics in IGA to boost computational efficiency and support complex real-world geometries. Indeed, various methods have been proposed in each of the two topics. In the area of local refinement, T-splines \cite{ref:sederberg03, ref:sederberg04, ref:bazilevs10} and (truncated) hierarchical B-splines \cite{ref:forsey88, ref:kraft97, ref:vuong11, ref:giannelli12} are two of the mostly studied methods. On the other hand, in the area of supporting EVs, recent efforts have been devoted to recovering optimal approximation properties, such as geometrically continuous construction \cite{ref:kapl15}, manifold splines \cite{ref:majeed17, ref:bercovier17}, the degenerated patch (D-patch) method \cite{ref:reif97, ref:tnguyen16, ref:toshniwal17}, blended construction in 3D \cite{ref:wei18} and subdivision methods \cite{ref:xli19, ref:wei20}. 

Moreover, several methods have also been proposed to support both local refinement and EVs at the same time, such as isogeometric spline forests \cite{ref:scott14}, truncated hierarchical Catmull-Clark subdivision \cite{ref:wei15a, ref:wei15b}, and analysis-suitable unstructured T-splines (ASUT-splines) with applications to Kirchhoff-Love shells \cite{ref:casquero20}. Particularly, ASUT-splines combine analysis-suitable T-splines \cite{ref:scott12, ref:veiga12, ref:li14} with the local D-patch construction \cite{ref:toshniwal17}, where the D-patch method is adopted to define $C^1$-continuous splines around EVs while using regular $C^2$-continuous B-splines away from EVs. While ASUT-splines show superior performance in Kirchhoff-Love shell problems \cite{ref:casquero20}, the need to resolve T-mesh constraints generally leads to refinement propagation beyond the region of interest. This becomes even more cumbersome when EVs are involved.

Therefore, we propose an alternative to ASUT-splines in this paper, namely \emph{truncated hierarchical unstructured splines} (THU-splines), by combining truncated hierarchical B-splines (THB-splines) with the local D-patch construction. Unlike T-splines, hierarchical B-splines (HB-splines) have a definite termination condition to drive local refinement. Moreover, there is no need for HB-splines to resolve mesh constraints to guarantee desired properties for the underlying spline functions. These advantages make HB-splines be able to locally refine regions of interest flexibly without propagation. However, applying THB-like refinement to the local D-patch construction is challenging because various types of spline functions are defined in a local D-patch construction. These different splines need to be well organized such that all hierarchical splines defined on each element can be evaluated in a unified manner. Also due to this challenge, proving properties of the underlying hierarchical basis becomes technically much more involving. While we postpone most of the theoretical study in a forthcoming paper, we study how THU-splines retain an invariant geometric mapping in a constructive manner, in order to provide insights of how the truncation in THU-splines works. Moreover, we will use examples to show the numerical evidence of other claimed properties, such as refinability and partition of unity.

The paper is organized as follows. We first briefly review two closely related methods, that is, THB-splines in Section~\ref{sec:thb} and local D-patch construction in Section~\ref{sec:local_Dpatch}. In Section~\ref{sec:thu}, we introduce our proposed method THU-splines. Section~\ref{sec:example} presents several numerical examples to show flexibility, effectiveness and efficiency of THU-splines. In the end, we conclude the paper and comment on the future work in Section~\ref{sec:conclude}.

\section{THB-splines}
\label{sec:thb}

We briefly review THB-splines in this section. THB-splines \cite{ref:giannelli12} were proposed as an enhancement of hierarchical B-splines (HB-splines) \cite{ref:kraft97} to reduce support overlapping of B-splines from different hierarchical levels, which in turn leads to a sparse stiffness matrix with reduced band width. Moreover, THB-spline basis functions are piecewise polynomials that form a non-negative partition of unity.

In what follows, we explain the key steps of THB-spline construction in the unvariate setting, including initialization, selection, and truncation. While these key ideas are the same in other dimension cases, interested readers are referred to \cite{ref:giannelli12, ref:lyche18} for details.

\subsection{Initialization}

To start with, a set of $n^0$ univariate B-splines $\mathcal{B}^0=\{B_i^0\}_{i=1}^{n^0}$ of degree $p$ ($n^0,\,p\in\mathbb{N}^{+}$), defined on the knot vector $\Xi^0 = \{\xi_i^0\}_{i=1}^{n^0+p+1}$, is given and initialized as the initial set of THB-spline basis functions, that is,
\begin{equation}
\mathcal{H}^0 = \mathcal{B}^0 .
\end{equation}
The initial domain serves as the entire domain $\Omega$ of THB-splines, i.e., $\Omega\equiv \Omega^0 = (\xi_1^0,\xi_{n+p+1}^0)$. 

To facilitate explanation, we introduce a series of nested spaces spanned by the B-splines,
\begin{equation}
\mathcal{B}^0 \subset \mathcal{B}^1 \subset \cdots \subset \mathcal{B}^{\ell_{\max}},	
\label{eq:nested_bsp}
\end{equation}
where $\ell_{\max}$ is the predefined maximum level, and the set $\mathcal{B}^{\ell}= \{B_i^\ell\}_{i=1}^{n^\ell}$ ($1\leq\ell\leq \ell_{\max}$) is defined on the knot vector $\Xi^\ell = \{\xi_i^\ell\}_{i=1}^{n^\ell+p+1}$. Note that the degree $p$ is fixed at all levels. As will become clear shortly, only a small subset of $\mathcal{B}^\ell$ is used and the full set $\mathcal{B}^\ell$ is never generated in practice. $\Xi^\ell$ is obtained by bisecting all the nonzero knot spans of $\Xi^{\ell-1}$, so knot vectors are nested as well,
\begin{equation}
\Xi^0 \subset \Xi^1 \subset \cdots \subset \Xi^{\ell_{\max}}.
\end{equation}

A hierarchical construction involves a nested sequence of hierarchical domains,
\begin{equation}
\Omega \equiv \Omega^0 \supseteq \Omega^1 \supseteq \cdots \supseteq \Omega^{\ell_{\max}}.
\label{eq:thb_nested_domain}
\end{equation}
This nested sequence implies a series of local refinements and it drives where to perform local refinement. For instance, once a set of THB-splines is initialized, we use $\mathcal{H}^0$ to solve the problem of interest. According to a certain a posteriori error estimator, a subdomain of $\Omega^0$, which is in fact the domain at Level 1 (i.e., $\Omega^1$), is marked to perform local refinement. Once $\Omega^1$ is determined, THB-splines $\mathcal{H}^1$ can be constructed accordingly. 

THB-splines are constructed level by level, so we only need to study a generic two-level construction for selection and truncation: given the THB-splines $\mathcal{H}^{\ell}$ that has been constructed up to Level $\ell$ ($0\leq \ell <\ell_{\max}$) and $\Omega^{\ell+1}$ (i.e., the subdomain of $\Omega^\ell$ to be refined), we discuss how to build $\mathcal{H}^{\ell+1}$.

\subsection{Selection}
\label{sec:thb_select}

The selection mechanism, first proposed in \cite{ref:kraft97} and later generalized in \cite{ref:vuong11}, is aimed to ensure linear independence of resulting hierarchical basis functions. It is based on the relation between the support of a candidate B-spline and \emph{active} subdomains of $\Omega^\ell$ and $\Omega^{\ell+1}$. The selection at Level $\ell$ states that a B-spline $B_i^\ell$ is selected only if its support has a nonzero intersection with $\Omega_a^\ell:=\Omega^\ell\backslash \Omega^{\ell+1}$, the \emph{active} subdomain at Level $\ell$. All the selected Level-$\ell$ B-splines are collected in $\mathcal{B}_a^\ell$,
\begin{equation}
\mathcal{B}_a^\ell = \{ B_i^\ell \in \mathcal{B}^\ell :\, \mathrm{supp}B_i^\ell \cap \Omega_a^\ell \neq \varnothing\}.
\label{eq:select_h0}
\end{equation}
where $\mathrm{supp}B_i^\ell := (\xi_i^\ell, \xi_{i+p+1}^\ell)$ is the support of $B_i^{\ell}$. On the other hand, at Level $\ell+1$, a B-spline is selected only if its support is fully contained in $\Omega^{\ell+1}$. All the selected Level-($\ell+1$) B-splines form the set $\mathcal{B}_{a,0}^{\ell+1}$,
\begin{equation}
\mathcal{B}_{a,0}^{\ell+1} = \{ B_i^{\ell+1} \in \mathcal{B}^{\ell+1} :\, \mathrm{supp}B_i^{\ell+1} \subseteq \Omega^{\ell+1} \}.
\label{eq:select_h}
\end{equation}
Selected B-splines are called \emph{active}, and the complementary set of B-splines (i.e., $\mathcal{B}^\ell \backslash \mathcal{B}_a^\ell$ and $\mathcal{B}^{\ell+1} \backslash \mathcal{B}_{a,0}^{\ell+1}$) are \emph{passive}. A prefix ``H" will be added (i.e., H-active or H-passive) to emphasize the hierarchical refinement when the context is not evident; see Section~\ref{sec:thu}. The subscript ``0" in $\mathcal{B}_{a,0}^{\ell+1}$ indicates the initial selection at Level $\ell+1$ with the entire $\Omega^{\ell+1}$ being active. Certain functions in $\mathcal{B}_{a,0}^{\ell+1}$ may latter become passive when Level $\ell+2$ is constructed from Level $\ell+1$ and the active Level-($\ell+1$) subdomain changes from $\Omega^{\ell+1}$ to $\Omega_a^{\ell+1}\equiv \Omega^{\ell+1}\backslash \Omega^{\ell+2}$. $\mathcal{I}^\ell$, $\mathcal{A}^\ell$, $\mathcal{I}^{\ell+1}$ and $\mathcal{A}_0^{\ell+1}$ are introduced to denote the index sets of $\mathcal{B}^\ell$, $\mathcal{B}_a^\ell$, $\mathcal{B}^{\ell+1}$ and $\mathcal{B}_{a,0}^{\ell+1}$, respectively. 

\subsection{Truncation}

The truncation mechanism, which enables reduced overlapping and partition of unity, plays the key role in THB-splines. We need the \emph{refinability relationship} to explain the idea. Refinability states that a Level-$\ell$ B-spline $B_i^\ell$ can be expressed as a linear combination of certain Level-($\ell+1$) B-splines,
\begin{equation}
B_i^\ell = \sum_{j \in \mathcal{C}_i^{\ell+1} } h_{ij}^{\ell+1} B_j^{\ell+1} , \quad \text{with} \quad \mathcal{C}_i^{\ell+1} := \{ j\in \mathcal{I}^{\ell+1} :\, \mathrm{supp}B_j^{\ell+1} \subset \mathrm{supp} B_i^\ell \} ,
\label{eq:refinability}
\end{equation}
where $h_{ij}^{\ell+1}$ are coefficients obtained from the knot insertion algorithm \cite{ref:boehm80}, and $B_j^{\ell+1}$ ($j\in\mathcal{C}_i^{\ell+1}$) are called the \emph{children} (or H-children, children in the hierarchical refinement) of $B_i^\ell$. In fact, refinability is the reason for Eq. \eqref{eq:nested_bsp} to hold.

The truncation mechanism applies to $B_i^\ell$ when some of its children are active, and it discards the active children from Eq. \eqref{eq:refinability}. As a result, the truncated function of $B_i^\ell$ is defined as
\begin{equation}
\mathrm{trun}_H B_i^\ell = \sum_{j \in \mathcal{C}_i^{\ell+1} \backslash \mathcal{A}_0^{\ell+1} } h_{ij}^{\ell+1} B_j^{\ell+1},
\label{eq:trun}
\end{equation}
which is also referred to as the \emph{inter-level truncation} in this paper.
Equivalently speaking, the coefficient $h_{ij}^{\ell+1}$ in Eq. \eqref{eq:refinability} is set to be $0$ by truncation if $B_j^{\ell+1}$ is active. Note that when all the children are passive (i.e., $\mathcal{C}_i^{\ell+1} \cap \mathcal{A}_0^{\ell+1} =\varnothing$), we have $\mathcal{C}_i^{\ell+1} \backslash \mathcal{A}_0^{\ell+1} \equiv \mathcal{C}_i^{\ell+1}$ and thus $\mathrm{trun}_H B_i^\ell \equiv B_i^\ell$. In other words, Eq. \eqref{eq:trun} can represent a generic definition including both truncated and non-truncated B-splines.

In summary, after selection of $\mathcal{B}_a^\ell$ and $\mathcal{B}_{a,0}^{\ell+1}$, truncation is applied to certain B-splines in $\mathcal{B}_a^\ell$, leading to a set of truncated (active) B-splines at Level $\ell$,
\begin{equation}
\mathcal{T}_a^{\ell} = \{ \mathrm{trun}_H B_i^\ell : \, B_i^\ell \in \mathcal{B}_a^\ell  \} .
\label{eq:tal}
\end{equation}
The THB-spline basis functions up to Level $\ell+1$ are now obtained as
\begin{equation}
\mathcal{H}^{\ell+1} = \mathcal{H}^{\ell-1} \cup \mathcal{T}_a^\ell \cup \mathcal{B}_{a,0}^{\ell+1} ,
\label{eq:THB_basis}
\end{equation}
where note that $\mathcal{H}^{\ell-1}:=\varnothing$ when $\ell=0$.

In Eq. \eqref{eq:THB_basis}, the local refinement is realized by enriching basis functions and it solely depends on $\mathcal{B}_{a,0}^{\ell+1}$. $\mathcal{H}^{\ell-1}$ and $\mathcal{T}_a^\ell$ are merely adjustments of existing basis functions to accommodate $\mathcal{B}_{a,0}^{\ell+1}$. Therefore, local refinement is effective only when $\mathcal{B}_a^{\ell+1} \neq \varnothing$; otherwise we end up with $\mathcal{H}^{\ell+1} \equiv \mathcal{H}^{\ell}$. A non-empty $\mathcal{B}_{a,0}^{\ell+1}$ implies that $\Omega^{\ell+1}$ must be able to cover the support of at least one Level-($\ell+1$) B-spline; see Eq. \eqref{eq:select_h}. This is indeed the only requirement for local refinement of THB-splines. It has a definite termination condition and it indicates a highly localized refinement with no propagation beyond the region of interest. 

\section{Local D-patch Construction}
\label{sec:local_Dpatch}

THB-splines are built upon B-splines that can only take structured meshes as input. However, our goal is to perform truncated hierarchical refinement on \emph{unstructured} quad meshes that involve EVs. As refinability is the key requirement for building hierarchical splines, we choose the D-patch method that retains this property while producing $C^1$-continuous splines in irregular regions, which we call \emph{D-patch splines} in the following. D-patch splines have shown optimal convergence in solving the 2nd order PDEs (partial differential equations), near-optimal convergence in the 4th order problems \cite{ref:toshniwal17}, as well as superior performance in Kirchhoff-Love shells \cite{ref:casquero20}.

In contrast to the \emph{global} D-patch construction \cite{ref:reif97, ref:tnguyen16} that leads to reduced continuity everywhere and proliferation of degrees of freedom, a \emph{local} D-patch construction defines D-patch splines only in irregular regions while keeping B-splines elsewhere. Such a local construction has been explored with global refinement \cite{ref:toshniwal17} as well as local refinement via T-splines \cite{ref:casquero20}.

Here we review the local D-patch construction following~\cite{ref:casquero20}, including B\'{e}zier extraction on unstructured meshes, D-patch construction around EVs, smooth transition to join D-patches with regular B-splines, and global refinement. Interested readers are referred to related literature \cite{ref:reif97, ref:toshniwal17, ref:casquero20} for further details.

\subsection{B\'{e}zier extraction}
\label{sec:bezier_ext}

\begin{figure}[htb]
\centering
\includegraphics[width=.5 \textwidth]{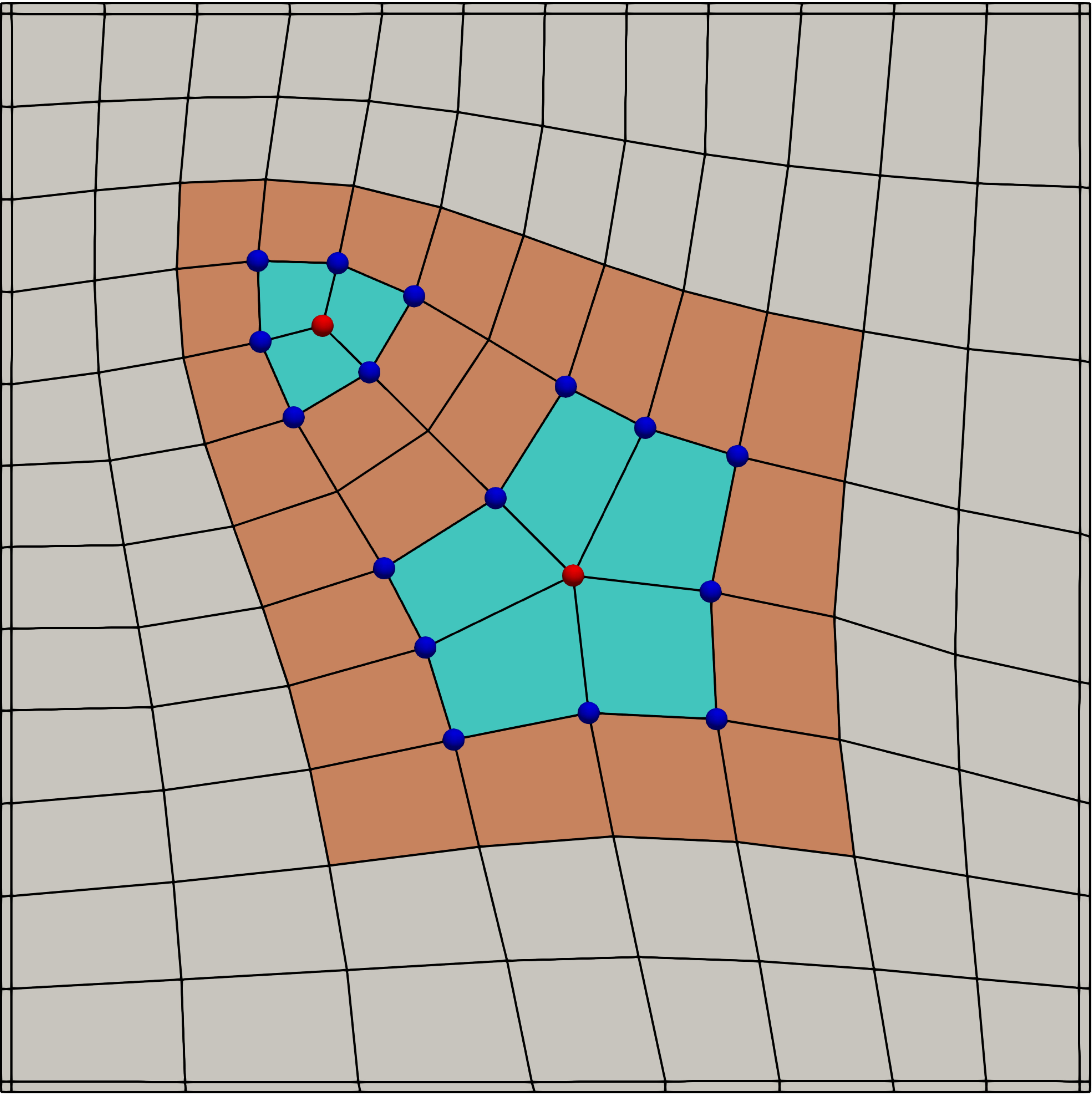}
\caption{An unstructured quad mesh taken as the input control mesh. EVs and interface vertices are marked as red and blue dots, respectively. Irregular elements are shaded blue and they are also transition elements in this mesh. Transition regular elements are shaded orange, whereas all the other elements are non-transition regular elements.}
\label{fig:unstruct_mesh}
\end{figure}

We first introduce several terminologies of unstructured meshes with the reference to Fig.~\ref{fig:unstruct_mesh}. An element is called an \emph{irregular element} if any of its vertices\footnote{We use ``vertex", ``edge" and ``face" to emphasize mesh connectivity (or topology), whereas using ``point" to carry the position or geometry information. ``Element" and ``face" are used interchangeably.} is an EV; otherwise it is regular. An element is a \emph{boundary element} if any of its vertices lies on the boundary; otherwise it is an interior element. The \emph{one-ring neighborhood} of a vertex includes all the elements sharing this vertex; recursively, the $n$-ring ($n\geq 2$) neighborhood is the ($n-1$)-ring neighborhood together with those sharing vertices with the ($n-1$)-ring neighborhood. The neighborhood of an element can be defined similarly. We adopt the following assumptions regarding the input mesh to simplify discussion.

\begin{itemize}
\item An irregular element can only have one EV;
\item A boundary vertex can only be shared by either one or two elements; and
\item An irregular element must not lie in the three-ring neighborhood of any boundary vertex.
\end{itemize}

Each vertex corresponds to a control point, and we call such points \emph{V-points} (i.e., vertex points).
Each edge is assigned with a non-negative real number called the knot span or knot interval, which is used to define element (parametric) domains as well as spline functions. For simplicity, we adopt a semi-uniform configuration of knot intervals where all edges have a unit knot interval except for those perpendicular to the boundary, which have a zero knot interval. As a result, edges around EVs all have the same knot intervals, which indeed is a necessary assumption for the D-patch construction \cite{ref:reif97}.

\begin{figure}[htb]
\centering
\includegraphics[width=.8 \textwidth]{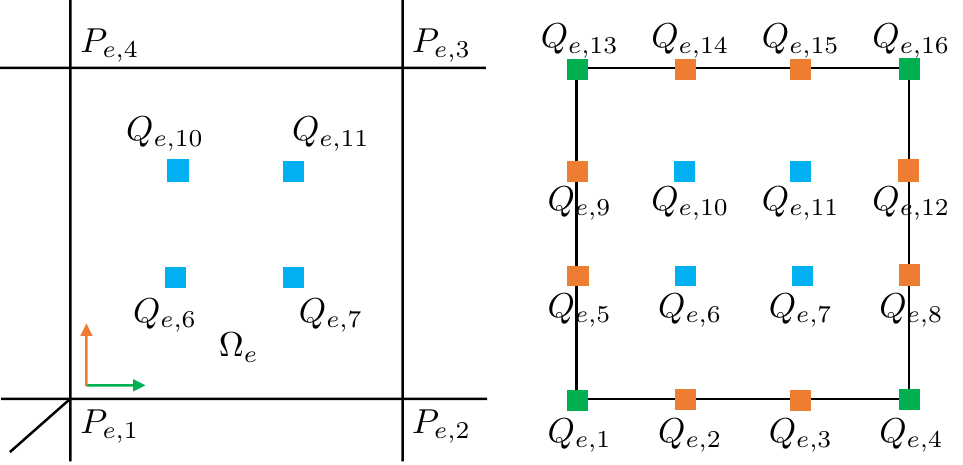}
\begin{tabular}{cc}
(a) & \hspace{-1.5cm} (b) \\
\includegraphics[width=.56 \textwidth]{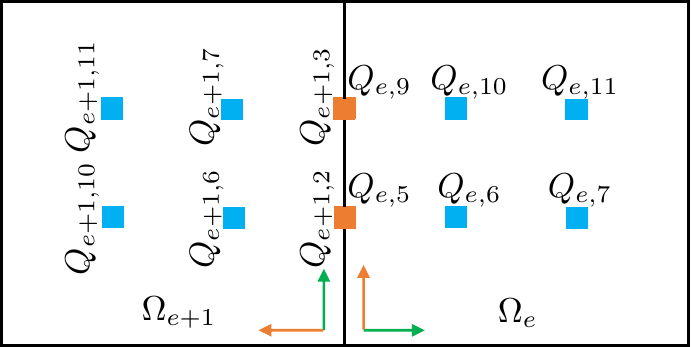} & \hspace{+1mm}
\includegraphics[width=.38 \textwidth]{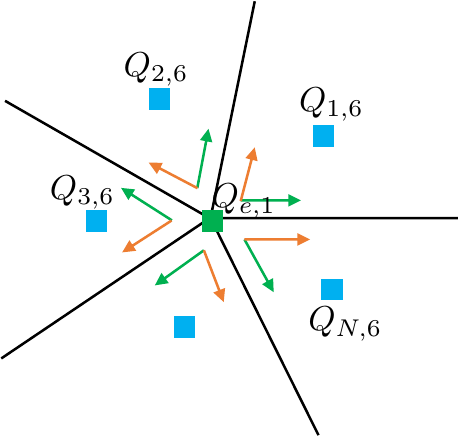} \\
(c) & (d) \\
\end{tabular}
\caption{Computing $4\times 4$ B\'{e}zier control points of $\Omega_e$. (a) The irregular element $\Omega_e$ and its four vertices $P_{e,i}$ ($i=1,\ldots,4$), where $P_{e,1}$ is an EV; (b) the 16 B\'{e}zier control points $Q_{e,j}$ ($j=1,\ldots,16$) of $\Omega_e$, where blue, orange and green squares are face, edge and vertex B\'{e}zier points, respectively; (c) computing an EB-point (edge B\'{e}zier point) as an average of neighboring FB-points (face B\'{e}zier points); and (d) computing a VB-point (vertex B\'{e}zier point) as an average of neighboring FB-points.}
\label{fig:bezier_ext}
\end{figure}

Note that given the third assumption regarding the input mesh, B\'ezier extraction is never needed for elements near the boundary. Without loss of generality, we next explain how to extract a B\'{e}zier representation for an irregular element $\Omega_e$ \cite{ref:scott13}.  Fig.~\ref{fig:bezier_ext}(a) shows as a generic irregular element $\Omega_e$, where $P_{e,i}$ ($i=1,\ldots,4$) are its 4 vertices, and $P_{e,1}$ is an EV. In the bicubic case, B\'{e}zier extraction is to find the $4\times 4$ B\'{e}zier control points $Q_{e,j}$ ($j=1,\ldots,16$) corresponding to $\Omega_e$. They can be divided into face, edge and vertex B\'{e}zier points, and they are called \emph{FB-points}, \emph{EB-points}, and \emph{VB-points}, respectively; see Fig. \ref{fig:bezier_ext}(b). The 4 FB-points are computed as convex combinations of the 4 element corners, that is,
\begin{equation}
\begin{bmatrix}
Q_{e,6} \\ Q_{e,7} \\ Q_{e,11} \\ Q_{e,10} \\
\end{bmatrix} =
\begin{bmatrix}
a & b & c & b \\
b & a & b & c \\
c & b & a & b \\
b & c & b & a \\
\end{bmatrix}
\begin{bmatrix}
P_{e,1} \\ P_{e,2} \\ P_{e,3} \\ P_{e,4} \\
\end{bmatrix} , \quad \text{with} \quad a=\frac{4}{9},\, b=\frac{2}{9},\, c=\frac{1}{9}.
\label{eq:face_pt}
\end{equation}
As shown in Fig. \ref{fig:bezier_ext}(c, d), EB- and VB-points are simply an average of their neighboring FB-points, for examples,
\begin{equation}
\begin{aligned}
Q_{e,9} &= Q_{e+1,3} = \frac{1}{2} Q_{e,10} + \frac{1}{2} Q_{e+1,7}, \\
Q_{e,1} &= \frac{1}{N} \sum_{i=1}^{N} Q_{i,6},
\end{aligned}
\label{eq:edge_vertex_pt}
\end{equation}
where $N\in\mathbb{N}^+$ is the number of faces sharing $P_{e,1}$, also called the valence of $P_{e,1}$.

Now let $\bm{Q}_e:=[Q_{e,1},\ldots,Q_{e,16}]^T$ denote the vector of B\'{e}zier points of $\Omega_e$, and $\bm{P}_v$ and $\bm{P}_f$ contain the V-points and FB-points of the one-ring neighborhood of $\Omega_e$, respectively. From Eqs. (\ref{eq:face_pt}, \ref{eq:edge_vertex_pt}), we observe that $\bm{Q}_e$ can be obtained through either
\begin{equation}
\bm{Q}_e = \bm{A}_e \bm{P}_{f},
\label{eq:bezier_pt_f}
\end{equation}
or
\begin{equation}
\bm{Q}_e = \bm{E}_e \bm{P}_v,
\label{eq:bezier_pt_v}
\end{equation}
where the matrix $\bm{A}_e$ contains the averaging coefficients in Eq. \eqref{eq:edge_vertex_pt}, and $\bm{E}_e$ is the extraction matrix with coefficients from both Eqs. (\ref{eq:face_pt}, \ref{eq:edge_vertex_pt}). When $N=4$, Eqs. (\ref{eq:bezier_pt_f}, \ref{eq:bezier_pt_v}) indicate B\'{e}zier extraction for a regular bicubic B-spline element. 

The B\'{e}zier representation of $\Omega_e$ is written as
\begin{equation}
\bm{S}(u,v)|_{\Omega_e} = \bm{Q}_e^T \bm{b}(u,v), \quad (u,v) \in [0,1]^2,
\label{eq:bezier_rep}
\end{equation}
where $\bm{b}(u,v):=[b_1(u,v),\ldots, b_{16}(u,v)]^T$ is the vector of bicubic Bernstein polynomials. $b_j(u,v)$ is expressed as,
\begin{equation}
b_j(u,v) = N_{(j-1)\%4}(u) N_{(j-1)/4}(v), \quad (u,v) \in [0,1]^2,
\end{equation}
where ``$\%$" and ``/" denote remainder and quotient operations, respectively; and
\begin{equation}
N_k(t) = \binom{3}{k} t^k (1-t)^{3-k}, \quad k=0,\ldots,3.
\end{equation}

Alternatively, the surface patch $\bm{S}(u,v)|_{\Omega_e}$ can also be expressed in terms of $\bm{P}_v$ and the associated spline functions $\bm{B}_v(u,v)$, 
\begin{equation}
\bm{S}(u,v)|_{\Omega_e} = \bm{P}_v^T \bm{B}_v(u,v).
\label{eq:spline_rep}
\end{equation}
On the other hand, Eqs. (\ref{eq:bezier_rep}, \ref{eq:spline_rep}) must be equivalent, leading to
\begin{equation}
\bm{P}_v^T \bm{B}_v(u,v) \equiv \bm{Q}_e^T \bm{b}(u,v) = \bm{P}_v^T \bm{E}_e^T \bm{b}(u,v),
\end{equation}
where we have used Eq. \eqref{eq:bezier_pt_v} to arrive at the last term.
According to the arbitrariness of $\bm{P}_v$, we conclude that
\begin{equation}
\bm{B}_v(u,v) = \bm{E}_e^T \bm{b}(u,v).
\label{eq:vspline_c0}
\end{equation}
Similarly, the surface patch can also be written with FB-points $\bm{P}_f$,
\begin{equation}
\bm{S}(u,v)|_{\Omega_e} = \bm{P}_f^T \bm{B}_f(u,v), \quad \text{with}\quad \bm{B}_f(u,v) = \bm{A}_e^T \bm{b}(u,v).
\label{eq:fspline_c0}
\end{equation}
The central idea of B\'{e}zier extraction is to express each spline function as a linear combination of Bernstein polynomials, which enables general spline functions to be easily integrated with existing finite element and isogeometric codes. While in regular elements, $\bm{B}_v$ and $\bm{B}_f$ are actually bicubic $C^2$ and $C^1$ B-splines, respectively, they are only $C^0$-continuous in irregular elements. The D-patch construction essentially adjusts $\bm{E}_e^T$ (or $\bm{A}_e^T$) to yield smooth spline functions, which we discuss in the following.

\subsection{D-patch method}
\label{subsec:dpatch}

In a local D-patch construction, the D-patch method is only adopted in regions around EVs, where we need to add necessary face-based control points and define their associated splines. Face-based control points are referred to as \emph{F-points} in the following. Splines associated with V-points and F-points are referred to as \emph{V-splines} and \emph{F-splines}, respectively.

\vspace{+2mm}\noindent\textbf{Enrichment via adding F-points}. 
We introduce \emph{enriched elements} to guide where to add F-points. Initially, all irregular elements are defined as enriched elements in the input mesh, whereas all the other elements are non-enriched elements. Enriched elements in refined meshes follow an ``inheritance" rule, which will be discussed in Section~\ref{sec:glb_refine}. 

\begin{figure}[htb]
\centering
\includegraphics[width=.5 \textwidth]{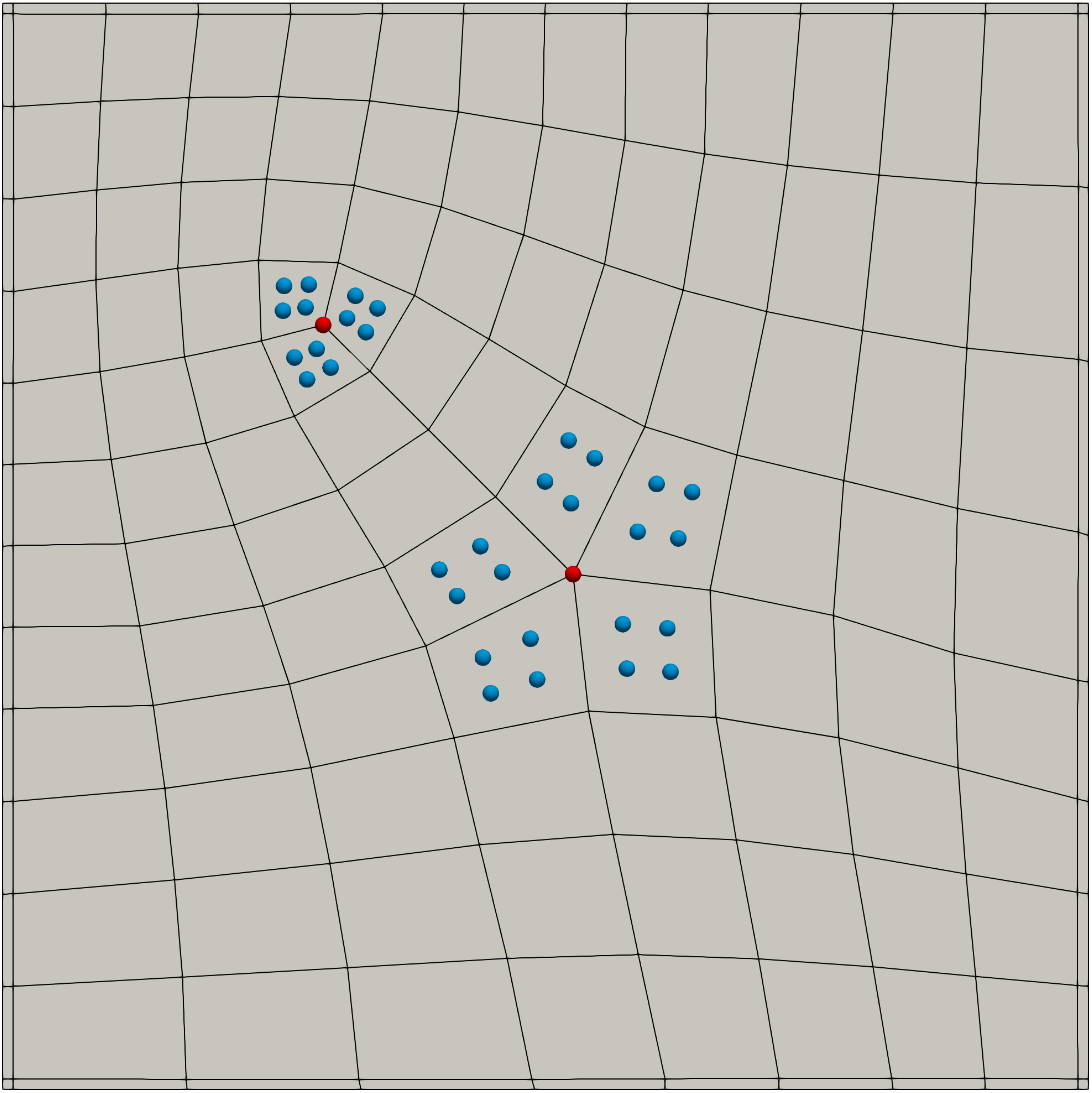}
\caption{The control mesh of a local D-patch representation corresponding to the input in Fig. \ref{fig:unstruct_mesh}, where blue dots are F-points, whereas red dots are EVs and they are passive.}
\label{fig:CM0}
\end{figure}

Four F-points are then added to each enriched element. In fact, F-points are equivalent to FB-points, and they are obtained according to Eq. \eqref{eq:face_pt} in the input mesh. As a result, the control mesh of a local D-patch representation consists of both V-points and F-points; see Fig.~\ref{fig:CM0}. However, while all the F-points of enriched elements are active in the sense that they are used in the local D-patch representation, not all V-points are active. A V-point is \emph{passive} if all its one-ring neighboring elements are enriched elements; otherwise it is \emph{active}. A prefix ``D" will be added to emphasize the context of local D-patch representation when necessary to distinguish those in hierarchical refinement. Clearly, all EVs are passive according to this definition. Essentially, passive V-points are replaced by their neighboring F-points. Note that such a control mesh corresponds to the analysis space discussed in \cite{ref:toshniwal17}.

\begin{figure}[htb]
\centering
\includegraphics[width=.7 \textwidth]{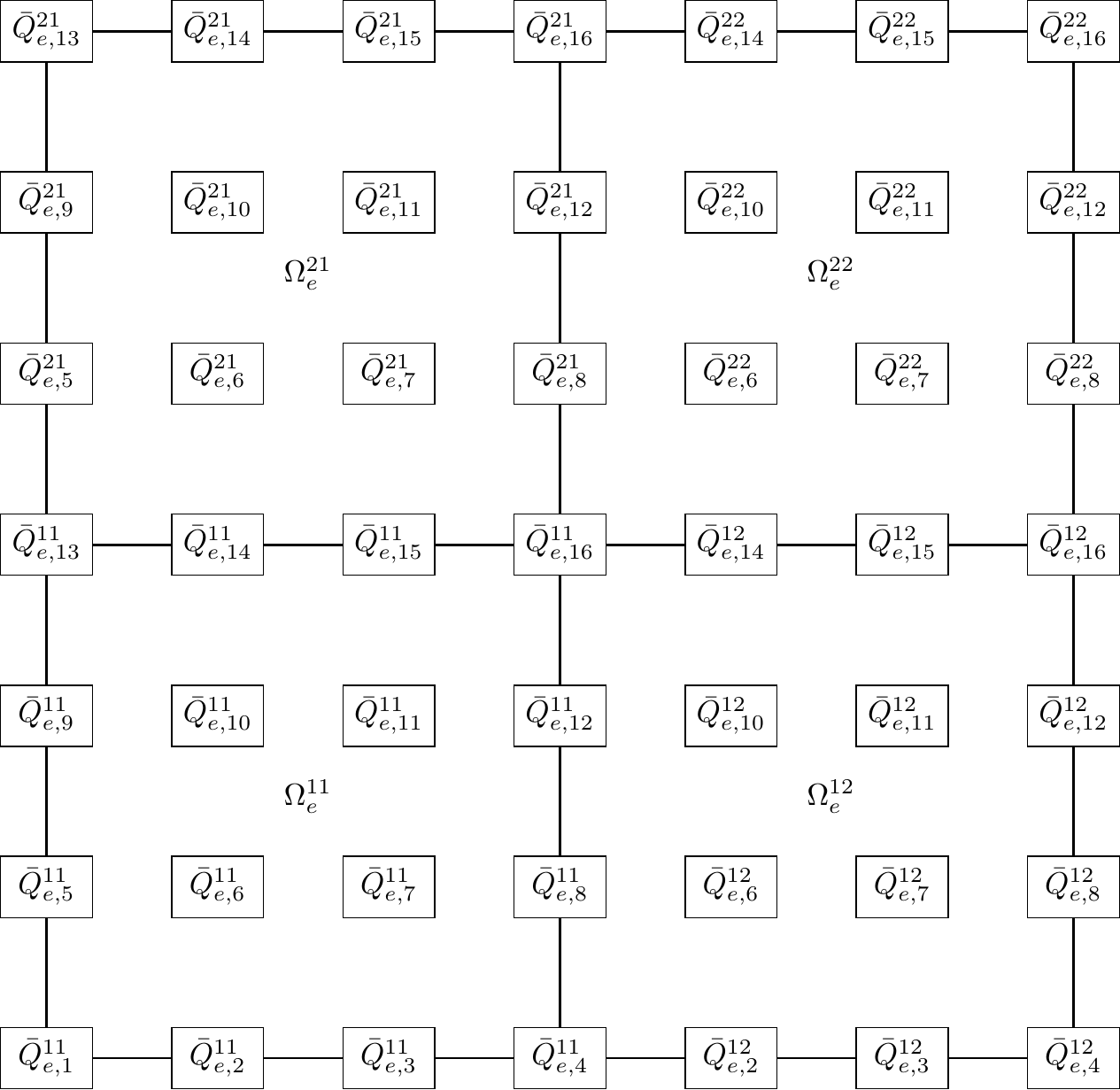}
\caption{$2\times 2$ split of the B\'{e}zier element in Fig. \ref{fig:bezier_ext}(b).}
\label{fig:bezier_split}
\end{figure}

\vspace{+2mm}\noindent\textbf{Definition of D-patch splines}.
We next briefly discuss the smooth D-patch splines around an EV that is shared by $N$ elements $\Omega_e$, $e=1,2,\ldots,N$. Interested readers are referred to \cite{ref:reif97, ref:toshniwal17, ref:casquero20} for details. D-patch splines are developed based on the $C^0$ version in Eq. \eqref{eq:fspline_c0}. This essentially involves smoothing the underlying surface patches around the EV, so we first proceed from the perspective of B\'{e}zier control points. Starting from Eq. \eqref{eq:bezier_pt_f}, we refine the B\'{e}zier patch via a $2\times 2$ split,
\begin{equation}
\bar{\bm{Q}}_{e}^{ij} = \bm{D}_{ij} \bm{Q}_e,
\label{eq:bezier_split}
\end{equation}
where $\bm{D}_{ij}$ ($i,j=1,2$) is the split matrix obtained from the well-known de Casteljau's algorithm \cite{ref:boehm99}, and $\bar{\bm{Q}}_{e}^{ij}$ represents control points of the refined B\'{e}zier patch corresponding to the subelement $\Omega_e^{ij}$; see Fig.~\ref{fig:bezier_split}. In particular, $\bar{Q}_{1,1}^{11}=\bar{Q}_{2,1}^{11}=\cdots=\bar{Q}_{N,1}^{11}=:\bar{Q}_{EV}$.

In the D-patch method, a smooth representation is achieved by adjusting all $\bar{Q}_{e,k}^{11}$, $e\in\{1,\ldots,N\}$, $k\in\{1,2,\ldots,16\}\backslash \{11,12,15,16\}$, which are collected into a vector $\bar{\bm{Q}}_{s}^{11}$. This is done through a linear projection matrix $\bm{\Pi}$,
\begin{equation}
\tilde{\bm{Q}}_s^{11} = \bm{\Pi} \bar{\bm{Q}}_s^{11}.
\end{equation}
As a result, for every $e\in\{1,\ldots,N\}$, $\tilde{Q}_{e,1}^{11}=\tilde{Q}_{e,2}^{11}=\tilde{Q}_{e,5}^{11}=\tilde{Q}_{e,6}^{11}\equiv\bar{Q}_{EV}$, and the five points $\tilde{Q}_{e,3}^{11}$, $\tilde{Q}_{e,7}^{11}$, $\tilde{Q}_{e,9}^{11}$, $\tilde{Q}_{e,10}^{11}$ and $\bar{Q}_{EV}$ are coplanar; see indices in Fig.~\ref{fig:bezier_split}. The reminder of $\bar{\bm{Q}}_{e}^{ij}$ stays the same, and the resulting B\'{e}zier surface patches are smoothly joined together around the EV.
The smoothing matrix $\bm{\Pi}$ corresponds to the idempotent version of the analysis space \cite{ref:casquero20, ref:toshniwal17}. Interested readers are referred there for further details. Note that although $\bm{\Pi}$ contains negative coefficients and thus leads to slightly negative spline functions, it guarantees refinability, which is the fundamental property required in a hierarchical spline construction.

Now we discuss the F-splines defined on an irregular element $\Omega_e$, which are associated with the F-points of the one-ring neighborhood of $\Omega_e$. They are piecewise-defined on the four subelements ($\Omega_e^{ij}$, $i,j=1,2$) as
\begin{equation}
\tilde{\bm{B}}_f(u,v)|_{\Omega_e^{ij}} = \bm{A}_e^T \bm{D}_{ij}^T \bm{\Pi}_{e,ij}^T \bm{b}(u,v).
\label{eq:fspline_c1}
\end{equation}
where the entries of $\bm{\Pi}_{e,ij}$ come from $\bm{\Pi}$ with an proper arrangement for $\Omega_e^{ij}$, and in particular, $\bm{\Pi}_{e,22}$ is an identity matrix. Comparing Eqs. (\ref{eq:fspline_c0}, \ref{eq:fspline_c1}), we observe that by introducing matrices $\bm{D}_{ij}$ and $\bm{\Pi}_{e,ij}$, $C^1$ continuity is built into $\tilde{\bm{B}}_f$. Moreover, we have $\tilde{\bm{B}}_f \equiv \bm{B}_f$ outside of the irregular element $\Omega_e$ thanks to the $2\times 2$ split, or equivalently speaking, the presence of $\bm{D}_{ij}$. Mostly importantly, it has also been proved that $\tilde{\bm{B}}_f$ are refinable in the sense of Eq. \eqref{eq:refinability}; see \cite{ref:tnguyen16}. In the reminder of the paper, we omit ``$\sim$" over $\tilde{\bm{B}}_f$ to keep notations simple.

\subsection{Smooth transition}
\label{sec:transition}

Next, we introduce how to smoothly join enriched and non-enriched elements in a local D-patch representation. The following terms are needed to facilitate the discussion. An interior edge is called an \emph{interface edge} if it is shared by an enriched element and a non-enriched element. Vertices of interface edges are interface vertices. An element is then called a \emph{transition element} if it contains an interface vertex. Note that transition elements can be regular or irregular; see Fig. \ref{fig:unstruct_mesh}. Now let $\mathcal{V}$ and $\mathcal{F}$ denote the index sets of V-points and F-points, respectively, and $\mathcal{V}_a$ be the active subset of $\mathcal{V}$.

We now discuss the splines associated with the active V-points ($\mathcal{V}_{a}$), which include interface points ($\mathcal{V}_{t}$) and the active regular points ($\mathcal{V}_{r}$) that do not lie on the interface, that is, $\mathcal{V}_a=\mathcal{V}_{t} \cup \mathcal{V}_{r}$. Splines associated with $\mathcal{V}_{r}$ are simply bicubic $C^2$ B-splines, whereas those associated with $\mathcal{V}_{t}$ need \emph{truncation} to enable a smooth transition between enriched and non-enriched elements. It is related but not the same as that introduced in the context of THB-splines. It was first proposed in blended B-splines \cite{ref:wei18} to recover the optimal convergence property when using unstructured hexahedral meshes. To distinguish, we refer to the truncation in a local D-patch representation as the \emph{same-level truncation}. Recall that it is called the inter-level truncation in the context of THB-splines.

\begin{figure}[htb]
\centering
\begin{tabular}{cc}
\includegraphics[width=.47 \textwidth]{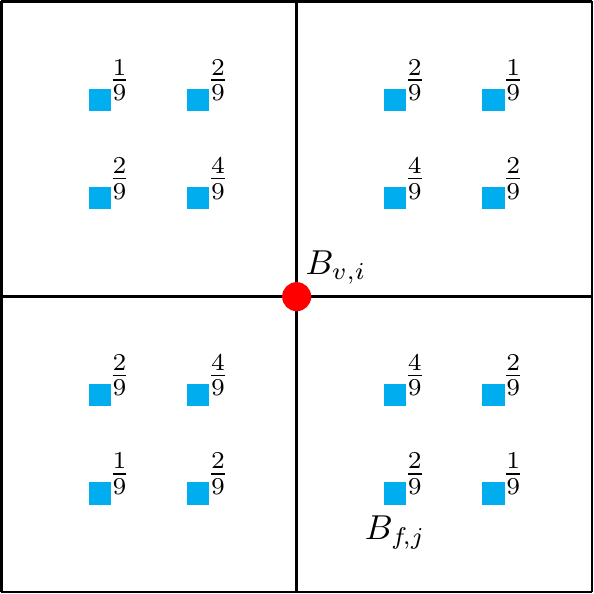} & \hspace{+1mm}
\includegraphics[width=.47 \textwidth]{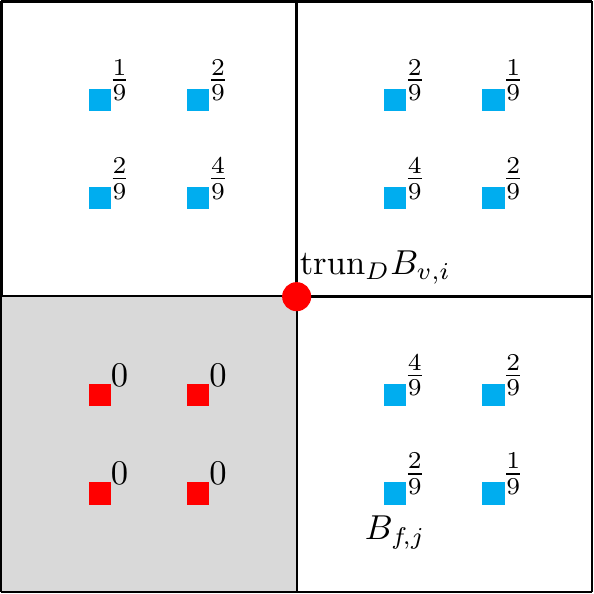} \\
(a) & (b) \\
\end{tabular}
\caption{The same-level truncation of $B_{v,i}$ that is associated with an interface point. (a) The V-spline $B_{v,i}$ is expressed in terms of one-ring neighboring F-splines $B_{f\!,j}$, where the numbers are the coefficients $d_{ij}$; and (b) truncation of $B_{v,i}$ by setting certain $d_{ij}=0$, where the gray element is enriched, and the red and blue squares are active and passive F-points, respectively.}
\label{fig:trun_u}
\end{figure}

For the V-spline $B_{v,i}$ associated with an interface point $P_i$ ($i\in\mathcal{V}_{t}$), we can express it in terms of its one-ring neighboring F-splines $B_{f\!,j}$,
\begin{equation}
B_{v,i} = \sum_{j\in\mathcal{F}_{i}} d_{ij} B_{f\!,j},
\label{eq:bv2f}
\end{equation}
where the coefficients $d_{ij}$ are obtained from Eq. \eqref{eq:face_pt} and they are illustrated in Fig. \ref{fig:trun_u}(a), $\mathcal{F}_{i}\subseteq \mathcal{F}$ is the index set of the F-points on the one-ring neighborhood of $P_i$, and $B_{f\!,j}$ are the children (or D-children, children in the local D-patch representation) of $B_{v,i}$. $B_{f\!,j}$ is defined by Eq. \eqref{eq:fspline_c1} if its F-point lies in an irregular element; otherwise it is defined by Eq. \eqref{eq:fspline_c0}. Note that not all $B_{f\!,j}$ are associated with the F-points of enriched elements. As $P_i$ is an interface point, there must exist at least one non-enriched element in its one-ring neighborhood. F-points on non-enriched elements are not part of the control mesh in a local D-patch representation, and thus they are passive. In contrast, F-points of enriched elements are active.

The same-level truncation for $B_{v,i}$ states that $d_{ij}$ is set to be 0 if $B_{f\!,j}$ is active. In other words, active children of $B_{v,i}$ are discarded,
\begin{equation}
\mathrm{trun}_D B_{v,i} = \sum_{j\in\mathcal{F}_{i} \backslash \mathcal{F}_a } d_{ij} B_{f\!,j},
\label{eq:trun_u}
\end{equation}
where $\mathcal{F}_a$ is the index set of active F-points. As a result of truncation, $\mathrm{trun}_D B_{v,i}$ is always a linear combination of regular $C^1$ B-splines.
Fig. \ref{fig:trun_u}(b) shows an example of $\mathrm{trun}_D B_{v,i}$, where the red squares are active F-points as the gray element is enriched whereas the blue squares are passive and only serve as an auxiliary purpose to define $\mathrm{trun}_D B_{v,i}$. Therefore, the splines associated with interface points are truncated B-splines defined according to Eq. \eqref{eq:trun_u}.  From an element-wise perspective, truncation has no influence outside of transition elements in the sense that $\mathrm{trun}_D B_{v,i}\equiv B_{v,i}$. 

On the other hand, given an active non-interface point $P_k$ ($k\in\mathcal{V}_{r}$), its associated spline $B_{v,k}$ is a $C^2$ B-spline and can be also expressed in the form of Eq. \eqref{eq:trun_u}, where, however, truncation is trivial as $\mathcal{F}_{k} \backslash \mathcal{F}_a = \mathcal{F}_{k}$, because we have $\mathcal{F}_{k} \cap \mathcal{F}_a = \varnothing$ by definition. In other words, $\mathrm{trun}_D B_{v,k}\equiv B_{v,k}$ holds for $k\in\mathcal{V}_{r}$.

To this end, all spline functions of a local D-patch construction have been discussed, which form a partition of unity and are linearly independent \cite{ref:casquero20}. The local D-patch representation is written as
\begin{equation}
S = \sum_{i\in\mathcal{F}_a} Q_i \, B_{f\!,i} + \sum_{i\in\mathcal{V}_{t}} P_i \, \mathrm{trun}_D B_{v,i} + \sum_{i\in\mathcal{V}_{r}} P_i \, B_{v,i} = \sum_{i\in\mathcal{F}_a} Q_i \, B_{f\!,i} + \sum_{i\in\mathcal{V}_a} P_i \, \mathrm{trun}_D B_{v,i},
\label{eq:surf_dpatch}
\end{equation}
where we recall that $\mathcal{V}_a = \mathcal{V}_r \cup \mathcal{V}_t$, and $Q_i$ and $P_i$ are (active) F-points and V-points, respectively.
The corresponding matrix form is given as,
\begin{equation}
\bm{S} =
\begin{bmatrix}
\bm{Q}^T & \bm{P}^T
\end{bmatrix} 
\begin{bmatrix}
\bm{B}_f \\ \bm{B}_v
\end{bmatrix},
\label{eq:dpatch_rep_glb}
\end{equation}
where $\bm{B}_v$ and $\bm{B}_f$ are vectors of V-splines and F-splines, respectively, and $\bm{B}_v$ includes both regular $C^2$ B-splines ($\mathcal{V}_r$) and truncated B-splines ($\mathcal{V}_t$).

\subsection{Global refinement}
\label{sec:glb_refine}

Finally, we discuss how to perform global $h$-refinement for a local D-patch representation. In the absence of EVs, $h$-refinement of a B-spline patch is straightforward through the knot insertion algorithm, where the underlying geometric mapping stays the same during refinement. We call such a refinement the \emph{consistent refinement}. However, consistent refinement is usually not available when EVs are involved. Currently, only two methods\footnote{We only focus on methods based on quadrilaterals.} have this property: the D-patch method and subdivision. In fact, subdivision implies a large family of methods that feature an infinite piecewise definition of splines in irregular regions. Interested readers are referred to \cite{ref:xli19, ref:wei20} for some of the latest developments that aim to recover optimal convergence rates.

Consistent global refinement of a local D-patch representation consists of a topological step and a geometric step. The topological step deals with mesh connectivity as well as inheritance of enriched elements, whereas the geometric step computes control point coordinates of the refined mesh. Elementwise definitions of basis functions are discussed as well after the topological step.

\vspace{+2mm}\noindent\textbf{Topological step}.
In the topological step, every element with a nonzero (parametric) measure is equally subdivided into 4 child elements. If a zero-measure element has a pair of opposite edges with nonzero knot intervals, it is subdivided into 2 child elements at the midpoints of such edges. Recall that irregular elements in the input mesh are marked as enriched elements. Subsequently in a refined mesh, elements in the ($2^k+1$)-ring neighborhood of each EV are marked as enriched elements, where $k\geq 1$ is number of global refinements. We observe the following: (1) regular elements in a refined mesh can be enriched elements; (2) child elements of a enriched element are always enriched elements; and (3) for a transition element that is not enriched, some of its child elements are also enriched, which are in fact the ($2^k+1$)$^{\text{th}}$-ring neighborhood of an EV. Such an inheritance pattern of enriched elements ensures that the irregular region in the initial local D-patch representation is maintained throughout refinements, which is the key to guaranteeing refinability~\cite{ref:toshniwal17}.

\vspace{+2mm}\noindent\textbf{Elementwise definitions of basis functions}.
Once topological step is finished, basis functions of a refined mesh can be defined accordingly. Their definitions vary among different element types: regular elements that are not transition or enriched (type-R), non-transition enriched elements (type-E), and transition elements (type-T). A type-R element only involves 16 $C^2$ B-splines. $C^1$-continuous splines are defined on a type-E element, which are D-patch splines if the element is irregular and 16 $C^1$ B-splines otherwise. On the other hand, a type-T element (no matter if the element is enriched or not) may involve a mix of all kinds of splines, such as V-splines (including $C^2$ B-splines and truncated B-splines) and F-splines (including $C^1$ B-splines and $C^1$ D-patch splines). As V-splines ($\bm{B}_v$) can be expressed as linear combinations of F-splines ($\bm{B}_f$), we adopt the following matrix form to compactly organize all these splines on a transition element,
\begin{equation}
\begin{bmatrix}
\bm{B}_v \\ \bm{B}_{f\!,a}
\end{bmatrix} =
\begin{bmatrix}
\bm{0} & \bm{T} \\
\bm{I} & \bm{0} \\
\end{bmatrix}
\begin{bmatrix}
\bm{B}_{f\!,a} \\ \bm{B}_{f\!,p} \\
\end{bmatrix} , \quad \text{with} \quad
\bm{B}_f = 
\begin{bmatrix}
\bm{B}_{f\!,a} \\ \bm{B}_{f\!,p} \\
\end{bmatrix},
\label{eq:trun_trans}
\end{equation}
where coefficients of the matrix $\bm{T}$ come from Eq. \eqref{eq:trun_u}, $\bm{I}$ is an identity matrix, and subscripts ``$a$" and ``$p$" again denote active and passive splines, respectively. 

\vspace{+2mm}\noindent\textbf{Geometric step}.
Positions of control points in a refined mesh are determined in the geometric step. Overall, all the coordinates are determined through the knot insertion algorithm. As has been discussed in \cite{ref:wei18}, this simple and unified computation is enabled by the (same-level) truncation. 

The computation of point coordinates has two cases: B-spline refinement for non-enriched elements, and B\'{e}zier refinement for enriched elements. First, for a surface patch restricted to a non-enriched element (no matter if the element is transition or not), it is simply a B-spline patch when we only focus on its V-points. All the V-points of refined patches are obtained through the knot insertion algorithm applied to such a B-spline patch. Moreover, when the non-enriched element is a transition element, recall that some of its child elements are enriched elements, whose F-points are added according to Eq. \eqref{eq:face_pt}.

Second, for a surface patch restricted to an enriched element, the corresponding local D-patch representation is first converted to its B\'{e}zier representation. Refinement is then applied to this B\'{e}zier patch using de Casteljau's algorithm. In each of the resulting child B\'{e}zier patches, the 4 interior B\'{e}zier control points are added to the corresponding element (which must be an enriched element by definition) as the F-points.

To this end, both basis functions and control points are defined for a refined mesh. Its local D-patch representation has the same form as Eq. \eqref{eq:dpatch_rep_glb}.

\section{THU-Splines}
\label{sec:thu}

We introduce the proposed method in this section, namely \emph{truncated hierarchical unstructured splines} (THU-splines). Following the construction of THB-splines, we discuss the three steps in constructing THU-splines: initialization, selection, and truncation. Among them, truncation becomes particularly complex due to the presence of two kinds of truncation: the same-level truncation and the inter-level truncation. 

\subsection{Initialization}

THU-splines are initialized with the local D-patch construction on an input unstructured quad mesh $\mathcal{M}^0$. For the ease of discussion, we introduce a series of consecutively refined meshes, $\{\mathcal{M}^{\ell}\}_{\ell=0}^{\ell_{\max}}$, where $\ell_{\max} \in\mathbb{N}^+$ is the predefined maximum number of refinements, and $\mathcal{M}^{\ell}$ is a global refinement of $\mathcal{M}^{\ell-1}$ ($1\leq\ell\leq \ell_{\max}$) according to Section~\ref{sec:glb_refine}. Note that global refinement is never actually performed and these globally refined meshes are introduced only to simplify explanation. Correspondingly, a local D-patch basis $\mathcal{B}^{\ell}$, which only consists of D-active splines, is constructed on each $\mathcal{M}^{\ell}$. The initial THU-spline basis $\mathcal{H}^0$ is given by $\mathcal{B}^0$. Moreover, with a closely related proof (Proposition 4.4, \cite{ref:toshniwal17}), we conjecture that refinability holds for $\{\mathcal{B}^\ell\}$, for which we will show numerical evidence in Section~\ref{sec:example}.

Let $\mathcal{F}^{\ell}$ and $\mathcal{V}_a^{\ell}$ denote the index sets of F-points and active V-points of $\mathcal{M}^{\ell}$, respectively. Note that $\mathcal{F}^{\ell}=\mathcal{F}_{a}^\ell \cup \mathcal{F}_{p}^\ell$, where passive F-points ($\mathcal{F}_{p}^\ell$) are introduced as an auxiliary means to define truncated B-splines associated with interface points, and splines associated with $\mathcal{F}_{p}^\ell$ do not belong to $\mathcal{B}^{\ell}$.

\subsection{Selection}

The same as in Section~\ref{sec:thb}, we focus on constructing $\mathcal{H}^{\ell+1}$ from $\mathcal{H}^{\ell}$. Following Eq.~\eqref{eq:thb_nested_domain}, we assume that a series of nested hierarchical domains $\{\Omega^\ell\}_{\ell=0}^{\ell_{\max}}$ is given to guide local refinement. However, we need to define the ``$\supset$" operation in Eq. \eqref{eq:thb_nested_domain} for unstructured meshes. Let $\Omega^0:=\{\Omega_e^0\}_{e=1}^{n^0}$ be a collection of all $n^0$ quad elements in $\mathcal{M}^0$. $\Omega^\ell$ ($1\leq\ell\leq\ell_{\max}$) is defined similarly but does not necessarily contain all elements of $\mathcal{M}^\ell$. $\Omega^{\ell} \supset \Omega^{\ell+1}$ means that every element in $\Omega^{\ell+1}$ is a child element of a certain element in $\Omega^\ell$. In other words, $\Omega^{\ell+1}$ is obtained by locally refining certain elements in $\Omega^{\ell}$; after refinement, the parent elements of those in $\Omega^{\ell+1}$ are no longer used and thus they are marked as passive.  The active subdomain of $\Omega_a^{\ell}$ consists of all active elements at Level $\ell$.

\begin{figure}[htb]
\centering
\begin{tabular}{cc}
\includegraphics[width=.47 \textwidth]{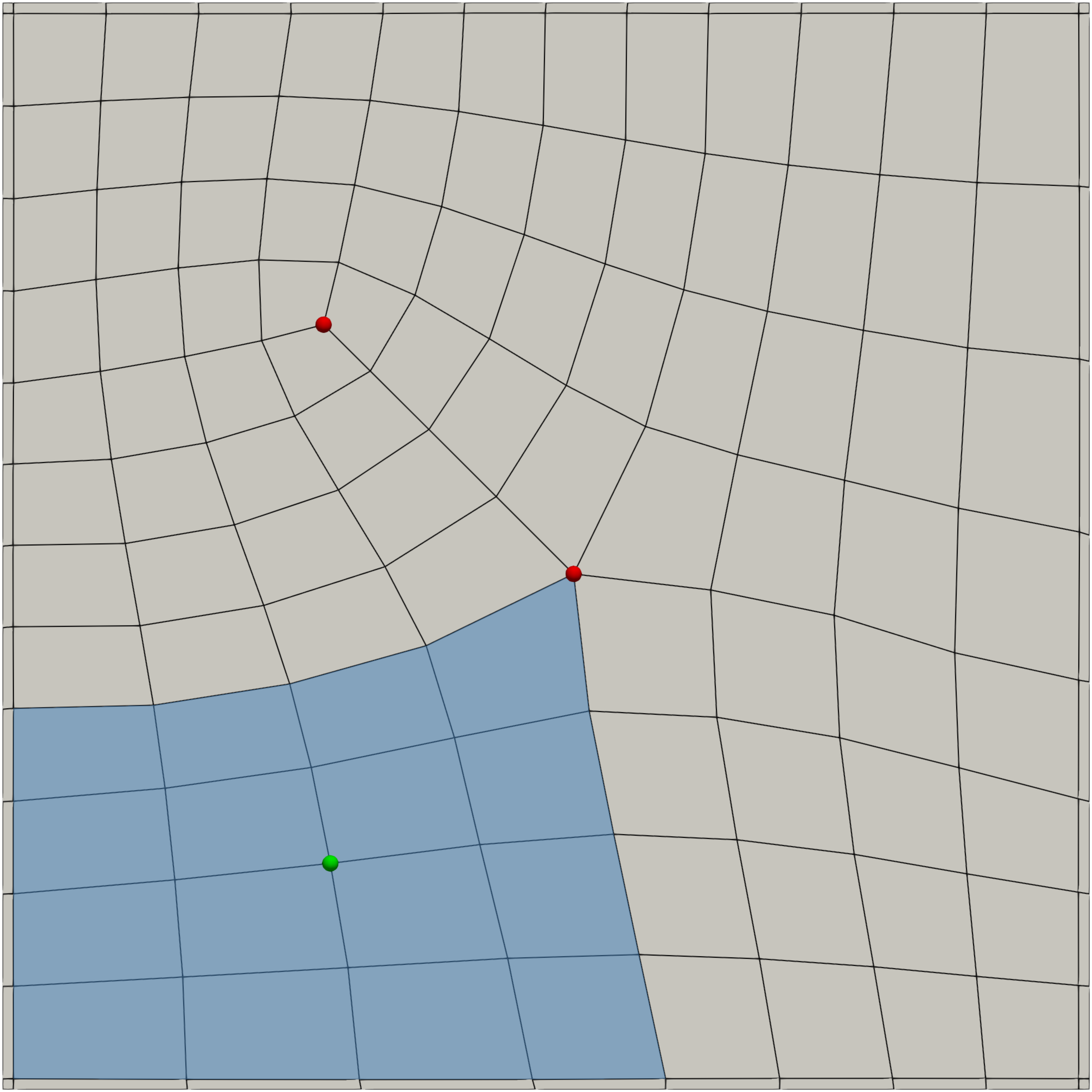} & \hspace{+1mm}
\includegraphics[width=.47 \textwidth]{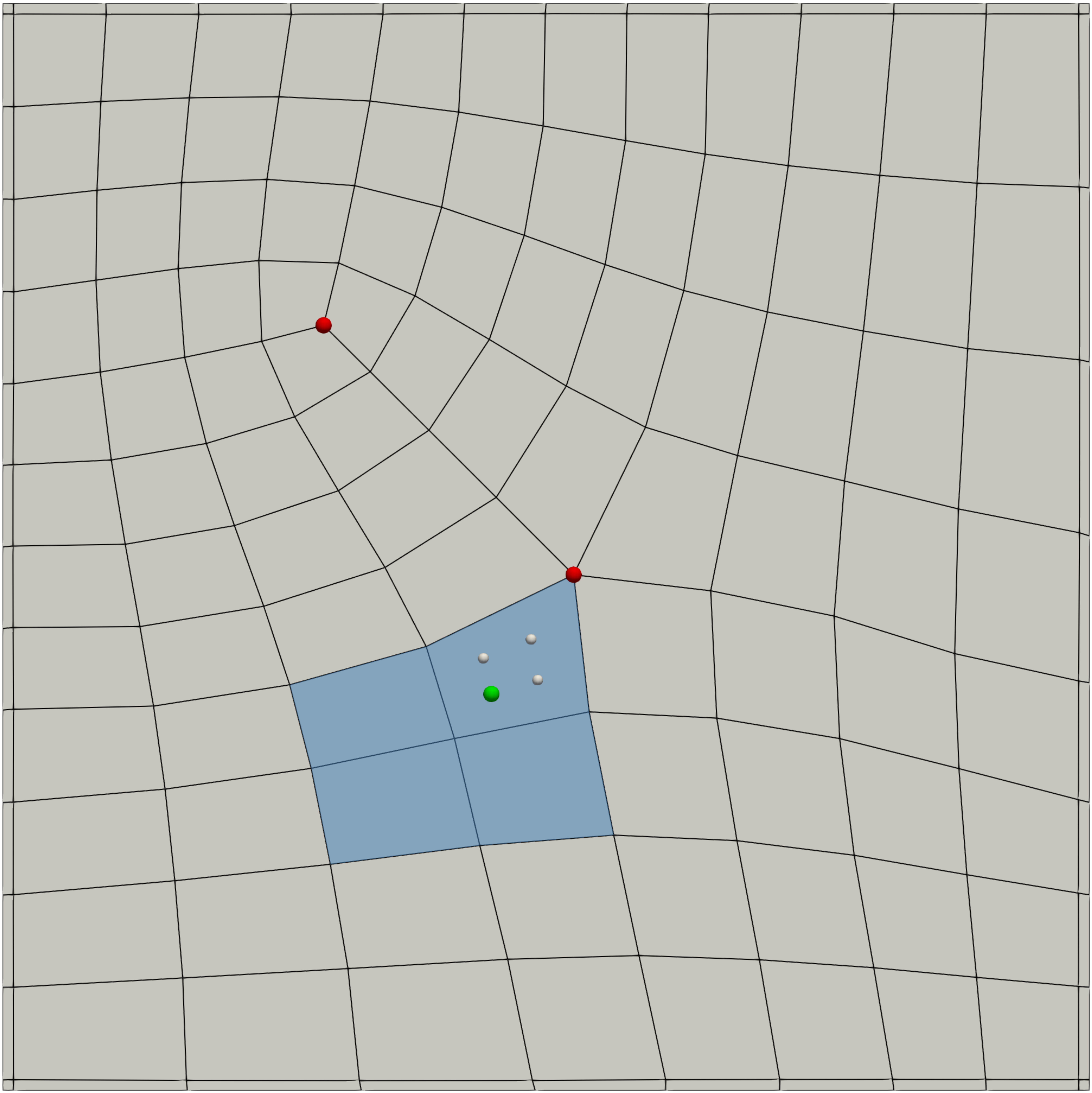} \\
(a) Case 1 & (b) Case 2 \\
\includegraphics[width=.47 \textwidth]{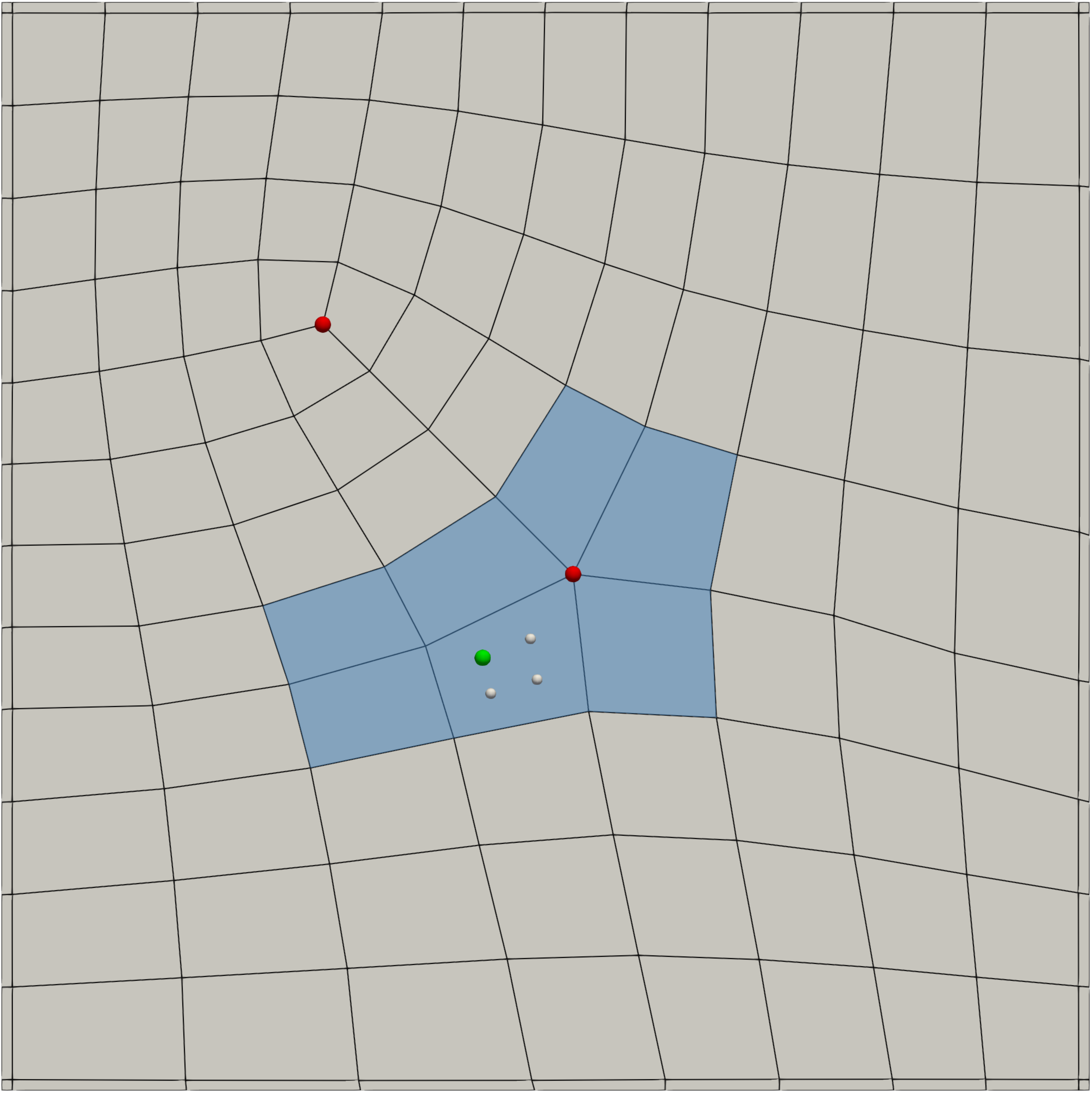} & \hspace{+1mm}
\includegraphics[width=.47 \textwidth]{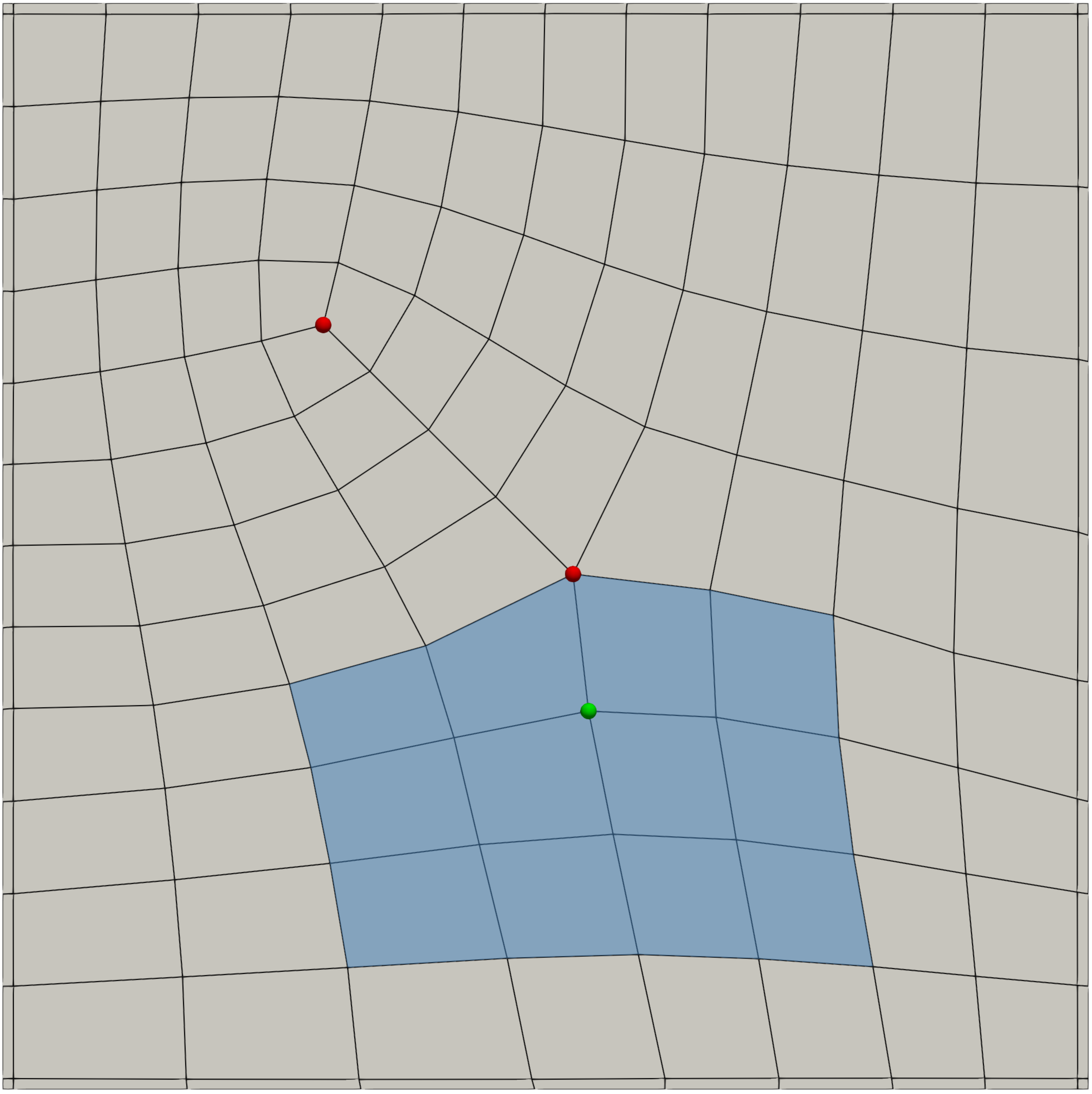} \\
(c) Case 3 & (b) Case 4 \\
\end{tabular}
\caption{Support of different spline functions, where red dots are extraordinary points, and blue shades are the supports of the green dots.}
\label{fig:support}
\end{figure}

It is left to define the support of every basis function to apply the selection mechanism. Given a generic basis function $B_i^\ell$ (whose control point is denoted as $P_i^\ell$), its support, $\mathrm{supp}B_i^\ell$, is defined as a collection of certain elements in $\mathcal{M}^\ell$. There are four cases depending on the type of $P_i^\ell$: (1) a V-point that is not on the interface, (2) an F-point on a regular element, (3) an F-point on an irregular element, and (4) an interface point; see Fig.~\ref{fig:support}.

In Case (1), $B_i^\ell$ is a bicubic $C^2$ B-spline, so $\mathrm{supp}B_i^\ell$ is the two-ring neighborhood of $P_i^\ell$. In Case (2), $B_i^\ell$ is a bicubic $C^1$ B-spline with doubly repeated knots, and $\mathrm{supp}B_i^\ell$ is the one-ring neighborhood of $P_i^\ell$. When $P_i^\ell$ is an F-point, its one-ring neighborhood means the one-ring neighborhood of its closest V-point. As a special case, the F-point $Q_{e,11}$ in Fig. \ref{fig:bezier_ext}(b), although being on an irregular element, also falls into Case (2), whereas the other three ($Q_{e,6}$, $Q_{e,7}$ and $Q_{e,10}$) belong to Case (3). 

Case (3) must involve a neighboring EV, and $B_i^\ell$ is defined as a D-patch spline. $\mathrm{supp}B_i^\ell$ consists of two parts. First, the same as Case (2), $B_i^\ell$ has support on the one-ring neighborhood of $P_i^\ell$. Moreover, $B_i^\ell$ is defined using the smoothing matrix such that it also has support on the one-ring neighborhood of the EV. The union of the two one-ring neighborhoods yields $\mathrm{supp}B_i^\ell$.

In Case (4), $B_i^\ell$ (more precisely, $\mathrm{trun}_D B_i^\ell$) is a truncated B-spline that is expressed as a linear combination of certain D-passive $C^1$ B-splines. Therefore, its support is the union of all the one-ring neighborhoods of these D-passive F-points.

\begin{figure}[htb]
\centering
\begin{tabular}{cc}
\includegraphics[width=.47 \textwidth]{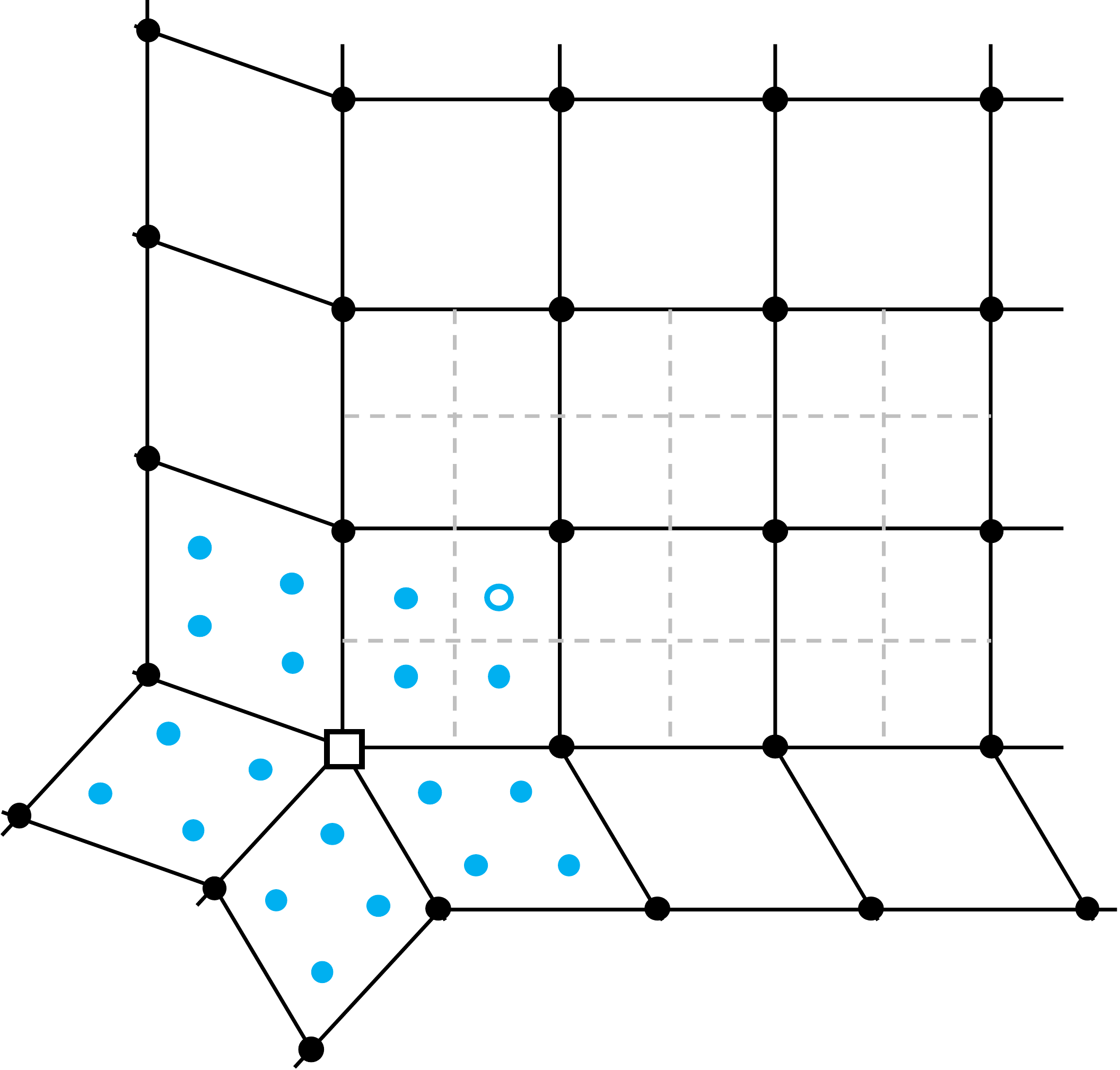} & \hspace{+1mm}
\includegraphics[width=.47 \textwidth]{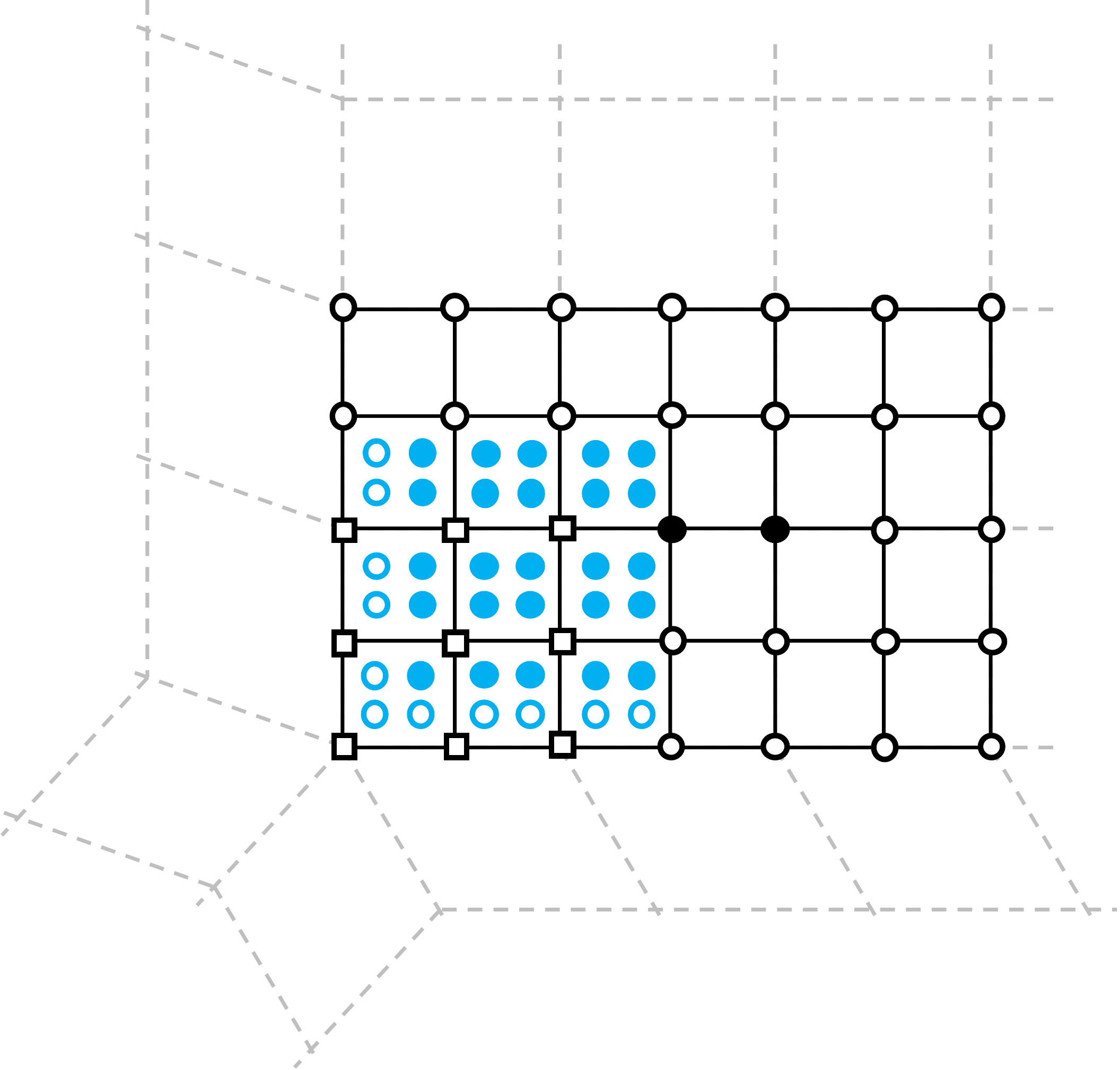} \\
(a) Level $\ell$ & (b) Level $\ell+1$ \\
\end{tabular}
\caption{Selection of spline functions at two consecutive levels. Dots and squares represent D-active and D-passive points, respectively. Filled and hollow dots are H-active and H-passive points, respectively. Black and blue dots (or squares) stand for V-points and F-points, respectively.}
\label{fig:select}
\end{figure}

With the hierarchical domains and the support of splines well defined, the selection mechanism follows the same as those in Eqs. (\ref{eq:select_h0}, \ref{eq:select_h}) to collect all the H-active splines as $\mathcal{B}_a^\ell$ and $\mathcal{B}_{a,0}^{\ell+1}$; see an example in Fig.~\ref{fig:select}. Note that during selection, only D-active splines at each level can be selected to be H-active or H-passive, whereas D-passive ones only serve as an auxiliary purpose to help define truncated B-splines and they never enter the selection procedure. In other words, D-passive functions are neither H-active nor H-passive.

Finally, recall that $\Omega^{\ell+1}$ needs to be large enough such that at least one Level-($\ell+1$) spline function can be added to $\mathcal{H}^{\ell+1}$. In the THU-spline construction, this condition can be easily guaranteed if the parent elements of those in $\Omega^{\ell+1}$ are a collection of one-ring neighborhoods of certain Level-$\ell$ points. In this manner, $\Omega^{\ell+1}$ at least contains the two-ring neighborhood of a Level-($\ell+1$) point. As the largest support in the above four cases is a two-ring neighborhood, a Level-($\ell+1$) spline can always be added.

\subsection{Truncation}

The inter-level truncation follows the same as in Eq. \eqref{eq:trun}. Certain H-active Level-$\ell$ functions are truncated. As a result, $\mathcal{B}_a^\ell$ is changed to $\mathcal{T}_a^\ell$, and $\mathcal{H}^{\ell+1}$ is constructed according to Eq. \eqref{eq:THB_basis}. Although the expressions stay the same as those in THB-splines, evaluation of THU-splines becomes much more involving because there are various types of splines rather than a single type as in THB-splines.

Without loss of generality, we study the THU-splines defined on a transition element $\Omega_e^{\ell+1}$, where we consider two levels (Levels $\ell$ and $\ell+1$) to simplify explanation. Following \cite{ref:bornermann13}, we evaluate all levels of H-active functions ($\bm{H}_{a}^{\ell}$ and $\bm{H}_{a}^{\ell+1}$) with the functions (both H-active and H-passive) at the highest level, and we write it in the matrix form, 
\begin{equation}
\bm{H}_a :=
\begin{bmatrix}
\bm{H}_{a}^\ell \\ \bm{H}_a^{\ell+1} \\
\end{bmatrix}=
\begin{bmatrix}
\bm{0} & \bm{C} \\
\bm{I} & \bm{0} \\
\end{bmatrix}
\begin{bmatrix}
\bm{H}_a^{\ell+1} \\ \bm{H}_p^{\ell+1} \\
\end{bmatrix},
\label{eq:thu_bf_0}
\end{equation}
where the coefficients of $\bm{C}$ come from Eq. \eqref{eq:trun}, and spline functions in $\bm{H}_{a}^{\ell}$ and $\bm{H}_{*}^{\ell+1}$ ($*\in\{a,p\}$) are from $\mathcal{B}_a^\ell$ and $\mathcal{B}^{\ell+1}$, respectively. Note that the inter-level truncation has been incorporated for $\bm{H}_a^{\ell}$ through the zero submatrix corresponding to $\bm{H}_a^{\ell+1}$. 

We further split $\bm{H}_a^\ell$ and $\bm{H}_*^{\ell+1}$ as V-splines ($\bm{H}_{v,a}^\ell$, $\bm{H}_{v,*}^{\ell+1}$) and F-splines ($\bm{H}_{f\!,a}^\ell$, $\bm{H}_{f\!,*}^{\ell+1}$). As a result, Eq. \eqref{eq:thu_bf_0} becomes
\begin{equation}
\begin{bmatrix}
\bm{H}_{v,a}^\ell \\ \bm{H}_{f\!,a}^\ell \\ \bm{H}_{v,a}^{\ell+1} \\ \bm{H}_{f\!,a}^{\ell+1} \\
\end{bmatrix}=
\begin{bmatrix}
\bm{0} & \bm{0} & \bm{C}_{11} & \bm{C}_{12} \\
\bm{0} & \bm{0} & \bm{C}_{21} & \bm{C}_{22} \\
\bm{I} & \bm{0} & \bm{0} & \bm{0} \\
\bm{0} & \bm{I} & \bm{0} & \bm{0} \\
\end{bmatrix}
\begin{bmatrix}
\bm{H}_{v,a}^{\ell+1} \\ \bm{H}_{f\!,a}^{\ell+1}  \\ \bm{H}_{v,p}^{\ell+1} \\ \bm{H}_{f\!,p}^{\ell+1} \\
\end{bmatrix},
\label{eq:thu_bf_1}
\end{equation}
where $\bm{C}_{21}=\bm{0}$ because an F-spline can only have F-splines as its children.
Moreover, $\bm{H}_{v,*}^{\ell+1}$ can be evaluated according to Eq. \eqref{eq:trun_trans}, that is, $\bm{H}_{v,*}^{\ell+1}=\bm{T}_{*}^{\ell+1} \bm{B}_{f\!,p}^{\ell+1}$, where $\bm{B}_{f\!,p}^{\ell+1}$ is a vector of D-passive F-splines that do not belong to $\mathcal{B}^{\ell+1}$. Substituting this into Eq. \eqref{eq:thu_bf_1}, we have
\begin{equation}
\bm{H}_a =
\begin{bmatrix}
\bm{H}_{v,a}^\ell \\ \bm{H}_{f\!,a}^\ell \\ \bm{H}_{v,a}^{\ell+1} \\ \bm{H}_{f\!,a}^{\ell+1} \\
\end{bmatrix}=
\begin{bmatrix}
\bm{0} & \bm{0} & \bm{C}_{11} & \bm{C}_{12} \\
\bm{0} & \bm{0} & \bm{0} & \bm{C}_{22} \\
\bm{I} & \bm{0} & \bm{0} & \bm{0} \\
\bm{0} & \bm{I} & \bm{0} & \bm{0} \\
\end{bmatrix}
\begin{bmatrix}
\bm{0} & \bm{0} & \bm{T}_{a}^{\ell+1} \\
\bm{I} & \bm{0} & \bm{0} \\
\bm{0} & \bm{0} & \bm{T}_{p}^{\ell+1} \\
\bm{0} & \bm{I} & \bm{0} \\
\end{bmatrix}
\begin{bmatrix}
\bm{H}_{f\!,a}^{\ell+1} \\ \bm{H}_{f\!,p}^{\ell+1} \\ \bm{B}_{f\!,p}^{\ell+1} \\
\end{bmatrix} =:
\bm{M}_f \bm{B}_f^{\ell+1},
\label{eq:thu_bf_2}
\end{equation}
where $\bm{B}_f^{\ell+1}$ contains the full set of F-splines (including both D-active ones $\bm{B}_{f\!,a}^{\ell+1}:= [\bm{H}_{f\!,a}^{\ell+1,T}, \bm{H}_{f\!,p}^{\ell+1,T}]^T$, and D-passive ones $\bm{B}_{f\!,p}^{\ell+1}$) defined on the transition element at Level $\ell+1$. Eq. \eqref{eq:thu_bf_2} implies that THU-splines $\bm{H}_a$ defined on a transition element can be evaluated in this unified manner as linear combinations of F-splines. The key is to obtain elementwise $\bm{M}_f$. With the two kinds of truncation, it can be shown that each column sum of $\bm{M}_f$ equals to one. Therefore, THU-splines $\bm{H}_a$ form a partition of unity. While we postpone the related proofs to a forthcoming paper, we will show numerical evidence in Section~\ref{sec:example}. 

When $\Omega_e^{\ell+1}$ is a non-transition enriched element (i.e., type-E), only F-splines $\bm{H}_{f\!,*}^{\ell+1}$ are involved. In other words, we have $\bm{H}_{v,*}^{\ell+1}=\bm{0}$ in Eq. \eqref{eq:thu_bf_1}. In contrast, when $\Omega_e^{\ell+1}$ is a type-R element, only $C^2$ B-splines $\bm{H}_{v,*}^{\ell+1}$ are involved and they not truncated. Moreover, $\bm{H}_{f\!,*}^{\ell+1}=\bm{0}$ as no F-splines are involved. In fact, this latter case coincides with standard THB-splines.

\subsection{Consistent refinement with THU-splines}
We discuss how THU-splines can achieve a consistent refinement in this section. We approach it in a constructive manner to show how the two kinds of truncation work together to retain this property.

We assume that only two levels (i.e., $\ell$ and $\ell+1$) are involved to provide insights of construction while keeping notations simple. Globally refined meshes ($\mathcal{M}^\ell$ and $\mathcal{M}^{\ell+1}$) and and their bases ($\mathcal{B}^\ell$ and $\mathcal{B}^{\ell+1}$) are in place to help discussion. Also recall the following notations for index sets at Level $L\in\{\ell,\ell+1\}$:
\begin{itemize}
\item $\mathcal{F}_a^L$ and $\mathcal{F}_p^L$: D-active and D-passive F-points/splines, respectively;
\item $\mathcal{V}_a^L$: D-active V-points/splines; and
\item $\mathcal{A}^L$: H-active points/splines (which can be F- or V-points/splines).
\end{itemize}

In what follows, we start with Level $\ell+1$ alone and then add back Level $\ell$. Without loss of generality, we study a transition element at Level $\ell+1$, $\Omega_e^{\ell+1}$, which can be enriched or non-enriched. Its surface patch can be written as follows with a full set of Level-($\ell+1$) F-splines,
\begin{equation}
S_e^{\ell+1} = \sum_{i\in \mathcal{F}_a^{\ell+1}} B_{f\!,i}^{\ell+1} Q_i^{\ell+1} + \sum_{i\in \mathcal{F}_p^{\ell+1}} B_{f\!,i}^{\ell+1} Q_i^{\ell+1},
\label{eq:surf_f}
\end{equation}
which is essentially the same as Eq.~\eqref{eq:fspline_c0}. Next, we express $S_e^{\ell+1}$ in terms of D-active points/splines only.
A D-passive F-point is obtained through its neighboring D-active V-points at the same level, so for $i\in\mathcal{F}_p^{\ell+1}$, we have
\begin{equation}
Q_i^{\ell+1} = \sum_{j\in\mathcal{V}_a^{\ell+1}} d_{ji}^{\ell+1} P_j^{\ell+1}, \quad \sum_{j\in\mathcal{V}_a^{\ell+1}} d_{ji}^{\ell+1} = 1,
\label{eq:Q_p}
\end{equation}
which is in fact the same as Eq.~\eqref{eq:bezier_pt_f}, with $d_{ji}^{\ell+1}\in\{4/9,2/9,1/9,0\}$. In fact, there are only four of $d_{ji}^{\ell+1}$ that are nonzero. Substituting Eq. \eqref{eq:Q_p} to Eq. \eqref{eq:surf_f}, we have
\begin{equation}
\begin{aligned}
S_e^{\ell+1} &= \sum_{i\in \mathcal{F}_a^{\ell+1}} B_{f\!,i}^{\ell+1} Q_i^{\ell+1}  + \sum_{i\in \mathcal{F}_p^{\ell+1}} B_{f\!,i}^{\ell+1} \sum_{j\in\mathcal{V}_a^{\ell+1}} d_{ji}^{\ell+1} P_j^{\ell+1} \\
&=\sum_{i\in \mathcal{F}_a^{\ell+1}} B_{f\!,i}^{\ell+1} Q_i^{\ell+1}  + \sum_{j\in\mathcal{V}_a^{\ell+1}} P_j^{\ell+1} \sum_{i\in \mathcal{F}_p^{\ell+1}} d_{ji}^{\ell+1}  B_{f\!,i}^{\ell+1} \\
&=\sum_{i\in \mathcal{F}_a^{\ell+1}} B_{f\!,i}^{\ell+1} Q_i^{\ell+1}  + \sum_{j\in\mathcal{V}_a^{\ell+1}} P_j^{\ell+1} \mathrm{trun}_D B_{v,j}^{\ell+1},
\end{aligned}
\label{eq:surf_l1}
\end{equation}
where note that we have used the definition of the same level in the last equation; see Eq. \eqref{eq:trun_u}. Eq. \eqref{eq:surf_l1} is the local D-patch representation and it is indeed equivalent to Eq. \eqref{eq:surf_dpatch}.

Now we add to Eq. \eqref{eq:surf_l1} the refinement relation between Levels $\ell+1$ and $\ell$. A Level-($\ell+1$) D-active V-point $P_j^{\ell+1}$ ($j\in\mathcal{V}_a^{\ell+1}$)  is computed with the Level-$\ell$ D-active V-points,
\begin{equation}
P_j^{\ell+1} = \sum_{j\in\mathcal{V}_a^\ell} h_{ji}^{\ell+1} P_j^{\ell},\quad \sum_{j\in\mathcal{V}_a^\ell} h_{ji}^{\ell+1} =1,
\label{eq:P_a}
\end{equation}
where $h_{ji}^{\ell+1}$ are the same coefficients as Eq. \eqref{eq:refinability} and only a small number of them are nonzero.
On the other hand, a Level-($\ell+1$) D-active F-point is computed from Level-$\ell$ F-points (both D-active and D-passive), so for $i\in\mathcal{F}_a^{\ell+1}$, we have
\begin{equation}
Q_i^{\ell+1} = \sum_{j\in\mathcal{F}_a^{\ell}} h_{ji}^{\ell+1} Q_j^{\ell} + \sum_{j\in\mathcal{F}_p^{\ell}} h_{ji}^{\ell+1} Q_j^{\ell}, \quad \sum_{j\in\mathcal{F}_a^{\ell}\cup\mathcal{F}_p^{\ell}} h_{ji}^{\ell+1} = 1,
\label{eq:Q_a}
\end{equation}
which is merely applying the B-spline knot insertion to doubly repeated knot vectors. When $j\in\mathcal{F}_p^\ell$, $Q_j^\ell$ is computed according to Eq.~\eqref{eq:Q_p} (replacing $\ell+1$ with $\ell$). 

We further divide $\mathcal{F}_a^{L}$ and $\mathcal{V}_a^{L}$ ($L\in\{\ell,\ell+1\}$) into H-active and H-passive subsets, leading to $\mathcal{A}_{f}^L:=\mathcal{F}_a^L \cap \mathcal{A}^L$, $\mathcal{A}_{v}^L:=\mathcal{V}_a^L \cap \mathcal{A}^L$, $\mathcal{N}_{f}^L:=\mathcal{F}_a^L \backslash \mathcal{A}^L$, and $\mathcal{N}_{v}^L:=\mathcal{V}_a^L \backslash \mathcal{A}^L$. Next, starting from Eq.~\eqref{eq:surf_l1}, we have
\begin{equation}
\begin{aligned}
S_e^{\ell+1} &= \sum_{i\in \mathcal{F}_a^{\ell+1}} B_{f\!,i}^{\ell+1} Q_i^{\ell+1}  + \sum_{j\in\mathcal{V}_a^{\ell+1}} P_j^{\ell+1} \mathrm{trun}_D B_{v,j}^{\ell+1} \\
&= \sum_{i\in \mathcal{A}_f^{\ell+1}} B_{f\!,i}^{\ell+1} Q_i^{\ell+1} + \sum_{i\in \mathcal{N}_f^{\ell+1}} B_{f\!,i}^{\ell+1} Q_i^{\ell+1} \\
&+ \sum_{j\in\mathcal{A}_v^{\ell+1}} P_j^{\ell+1} \mathrm{trun}_D B_{v,j}^{\ell+1} + \sum_{j\in\mathcal{N}_v^{\ell+1}} P_j^{\ell+1} \mathrm{trun}_D B_{v,j}^{\ell+1} .
\end{aligned}
\label{eq:sap}
\end{equation}
We keep the terms corresponding to H-active subsets unchanged but further expand H-passive ones using Eqs.~(\ref{eq:P_a}, \ref{eq:Q_a}),
\begin{equation}
\begin{aligned}
S_{e,\mathcal{N}}^{\ell+1} &:= \sum_{i\in \mathcal{N}_f^{\ell+1}} B_{f\!,i}^{\ell+1} Q_i^{\ell+1} + \sum_{j\in\mathcal{N}_v^{\ell+1}} P_j^{\ell+1} \mathrm{trun}_D B_{v,j}^{\ell+1} \\
&= \sum_{i\in \mathcal{N}_f^{\ell+1}} B_{f\!,i}^{\ell+1} \left( \sum_{j\in\mathcal{F}_a^{\ell}} h_{ji}^{\ell+1} Q_j^{\ell} + \sum_{j\in\mathcal{F}_p^{\ell}} h_{ji}^{\ell+1} Q_j^{\ell} \right) + \sum_{i\in\mathcal{N}_v^{\ell+1}} \left( \sum_{j\in\mathcal{V}_a^\ell} h_{ji}^{\ell+1} P_j^{\ell} \right) \mathrm{trun}_D B_{v,i}^{\ell+1}  \\
&=\sum_{i\in \mathcal{N}_f^{\ell+1}} B_{f\!,i}^{\ell+1} \left( \sum_{j\in\mathcal{F}_a^{\ell}} h_{ji}^{\ell+1} Q_j^{\ell} + \sum_{j\in\mathcal{F}_p^{\ell}} h_{ji}^{\ell+1} \sum_{k\in\mathcal{V}_a^\ell} d_{kj}^{\ell+1} P_k^{\ell} \right) \\
&+ \sum_{i\in\mathcal{N}_v^{\ell+1}} \left( \sum_{j\in\mathcal{V}_a^\ell} h_{ji}^{\ell+1} P_j^{\ell} \right) \mathrm{trun}_D B_{v,i}^{\ell+1}.
\end{aligned}
\label{eq:sp}
\end{equation}

Notice that according to the selection mechanism of hierarchical splines, a Level-$\ell$ spline (F- or V-spline) is H-passive if and only if all its children are H-active; in other words, if any of its children is passive, then the spline is active. Therefore in the last equation in Eq.~\eqref{eq:sp}, we observe that when $i\in\mathcal{N}_f^{\ell+1}$ and $h_{ji}^{\ell+1}\neq 0$, $Q_j^\ell$ must be H-active, i.e., $j\in\mathcal{A}_f^\ell$. Similarly, when $i\in\mathcal{N}_f^{\ell+1}$, $h_{ji}^{\ell+1}\neq 0$ and $d_{kj}^{\ell+1}\neq 0$, we have $k\in\mathcal{A}_v^\ell$ for $P_k^\ell$; when $i\in\mathcal{N}_v^{\ell+1}$ and $h_{ji}^{\ell+1}\neq 0$, we must have $j\in\mathcal{A}_v^\ell$ for $P_j^\ell$. Then after proper rearrangement, Eq.~\eqref{eq:sp} becomes
\begin{equation}
\begin{aligned}
S_{e,\mathcal{N}}^{\ell+1} &=\sum_{j\in\mathcal{A}_f^{\ell}} Q_j^{\ell}  \sum_{i\in \mathcal{N}_f^{\ell+1}} h_{ji}^{\ell+1} B_{f\!,i}^{\ell+1}  + \sum_{k\in\mathcal{A}_v^\ell} P_k^{\ell}  \sum_{i\in \mathcal{N}_f^{\ell+1}} B_{f\!,i}^{\ell+1} \sum_{j\in\mathcal{F}_p^{\ell}} h_{ji}^{\ell+1} d_{kj}^{\ell+1}  \\
&+ \sum_{j\in\mathcal{A}_v^\ell} P_j^{\ell}  \sum_{i\in\mathcal{N}_v^{\ell+1}}  h_{ji}^{\ell+1} \mathrm{trun}_D B_{v,i}^{\ell+1}\\
&=\sum_{j\in\mathcal{A}_f^{\ell}} Q_j^{\ell}  \sum_{i\in \mathcal{N}_f^{\ell+1}} h_{ji}^{\ell+1} B_{f\!,i}^{\ell+1}  + \sum_{k\in\mathcal{A}_v^\ell} P_k^{\ell}  \sum_{i\in \mathcal{N}_f^{\ell+1}} B_{f\!,i}^{\ell+1} \sum_{j\in\mathcal{F}_p^{\ell}} h_{ji}^{\ell+1} d_{kj}^{\ell+1}  \\
&+ \sum_{j\in\mathcal{A}_v^\ell} P_j^{\ell}  \sum_{i\in\mathcal{N}_v^{\ell+1}}  h_{ji}^{\ell+1} \sum_{k\in\mathcal{F}_p^{\ell+1}} d_{ik}^{\ell+1} B_{f\!,k}^{\ell+1},
\end{aligned}
\label{eq:sp1}
\end{equation}
where we have used the definition of $\mathrm{trun}_D B_{v,i}^{\ell+1}$ in the last term.
Let $s_{ki}^{\ell+1}:=\sum_{j\in\mathcal{F}_p^{\ell}} h_{ji}^{\ell+1} d_{kj}^{\ell+1}$ and $t_{jk}^{\ell+1}:=\sum_{i\in\mathcal{N}_v^{\ell+1}} h_{ji}^{\ell+1} d_{ik}^{\ell+1}$. We have
\begin{equation}
\begin{aligned}
S_{e,\mathcal{N}}^{\ell+1} &=\sum_{j\in\mathcal{A}_f^{\ell}} Q_j^{\ell}  \sum_{i\in \mathcal{N}_f^{\ell+1}} h_{ji}^{\ell+1} B_{f\!,i}^{\ell+1} + \sum_{k\in\mathcal{A}_v^\ell} P_k^{\ell}  \sum_{i\in \mathcal{N}_f^{\ell+1}} s_{ki}^{\ell+1} B_{f\!,i}^{\ell+1} + \sum_{j\in\mathcal{A}_v^\ell} P_j^{\ell}  \sum_{k\in\mathcal{F}_p^{\ell+1}}  t_{jk}^{\ell+1} B_{f\!,k}^{\ell+1}  \\
&=\sum_{j\in\mathcal{A}_f^{\ell}} Q_j^{\ell}  \sum_{i\in \mathcal{N}_f^{\ell+1}} h_{ji}^{\ell+1} B_{f\!,i}^{\ell+1} + \sum_{j\in\mathcal{A}_v^\ell} P_j^{\ell}  \left( \sum_{i\in \mathcal{N}_f^{\ell+1}} s_{ji}^{\ell+1} B_{f\!,i}^{\ell+1} + \sum_{i\in\mathcal{F}_p^{\ell+1}}  t_{ji}^{\ell+1} B_{f\!,i}^{\ell+1} \right) \\
&=\sum_{j\in\mathcal{A}_f^{\ell}} Q_j^{\ell}  \mathrm{trun}_H B_{f\!,j}^{\ell}  + \sum_{j\in\mathcal{A}_v^\ell} P_j^{\ell}  \mathrm{trun}_{H\!D} B_{v,j}^{\ell},
\end{aligned}
\label{eq:sp2}
\end{equation}
with
\begin{equation}
\mathrm{trun}_H B_{f\!,j}^\ell = \sum_{i\in \mathcal{N}_f^{\ell+1}} h_{ji}^{\ell+1} B_{f\!,i}^{\ell+1},
\label{eq:trun_f}
\end{equation}
and
\begin{equation}
\mathrm{trun}_{H\!D} B_{v,j}^{\ell} := \sum_{i\in \mathcal{N}_f^{\ell+1}} s_{ji}^{\ell+1} B_{f\!,i}^{\ell+1} + \sum_{i\in\mathcal{F}_p^{\ell+1}}  t_{ji}^{\ell+1} B_{f\!,i}^{\ell+1}.
\label{eq:trun_hd}
\end{equation}
Eq.~\eqref{eq:trun_f} is merely a THB-spline that has a knot vector with doubly repeated knots, whereas Eq.~\eqref{eq:trun_hd} implies a mix of the same-level truncation and the inter-level truncation encoded in the coefficients $s_{ji}^{\ell+1}$ and $t_{ji}^{\ell+1}$. Moreover, Eqs.~(\ref{eq:trun_f}, \ref{eq:trun_hd}) correspond to $\bm{H}_{f\!,a}^\ell$ and $\bm{H}_{v,a}^\ell$ in Eq.~\eqref{eq:thu_bf_2}, respectively. Although truncated splines in Eqs.~(\ref{eq:trun_f}, \ref{eq:trun_hd}) have different forms, the key idea is essentially the same: A truncated spline is always expressed as a linear combination of certain passive splines.

Finally, with Eqs.~(\ref{eq:sap}, \ref{eq:sp2}), the surface patch $S_e^{\ell+1}$ of the element $\Omega_e^{\ell+1}$ can be written in terms of H-active splines only,
\begin{equation}
\begin{aligned}
S_{e}^{\ell+1} &= \sum_{i\in \mathcal{A}_f^{\ell+1}} B_{f\!,i}^{\ell+1} Q_i^{\ell+1} + \sum_{i\in\mathcal{A}_v^{\ell+1}} P_i^{\ell+1} \mathrm{trun}_D B_{v,i}^{\ell+1} \\
&+ \sum_{i\in\mathcal{A}_f^{\ell}} Q_i^{\ell}  \mathrm{trun}_H B_{f\!,i}^{\ell}  + \sum_{i\in\mathcal{A}_v^\ell} P_i^{\ell}  \mathrm{trun}_{H\!D} B_{v,i}^{\ell},
\end{aligned}
\label{eq:thu_rep}
\end{equation}
which is also the general form of a THU-spline representation (with two levels). When $\Omega_e^{\ell+1}$ is a Type-R (non-transition and non-enriched) element, only THB-splines are involved, i.e., $\mathcal{A}_f^L=\varnothing$, $\mathrm{trun}_D B_{v,i}^{\ell+1}= B_{v,i}^{\ell+1}$ and $\mathrm{trun}_{H\!D} B_{v,i}^{\ell}= B_{v,i}^{\ell}$. When $\Omega_e^{\ell+1}$ is a Type-E (non-transition and enriched) element, only F-splines are involved, i.e., $\mathcal{A}_v^L=\varnothing$.

\section{Numerical Examples}
\label{sec:example}

\begin{figure}[htb]
\centering
\begin{tabular}{cc}
\includegraphics[width=.4 \textwidth]{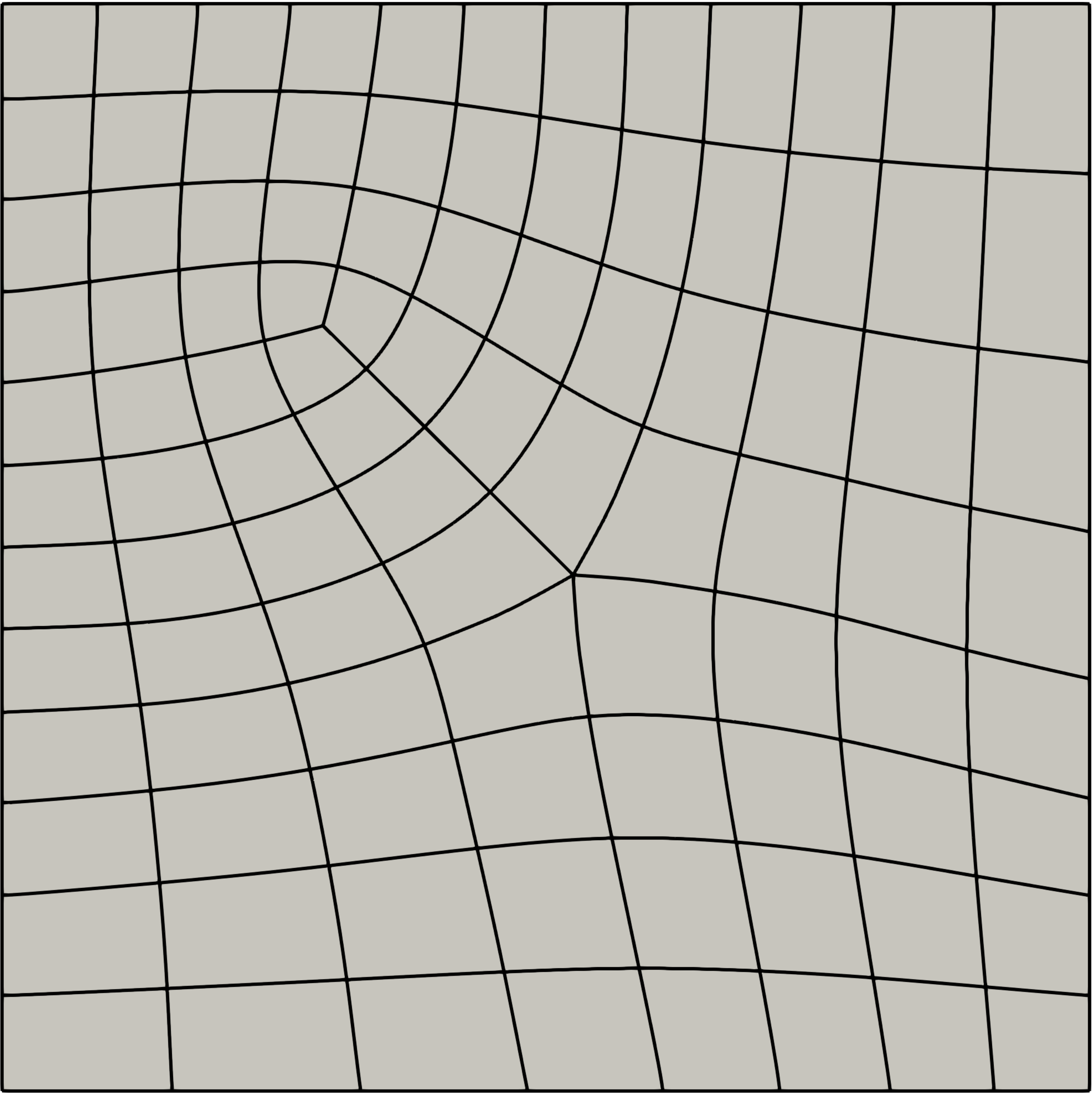} & \hspace{+1mm}
\includegraphics[width=.4 \textwidth]{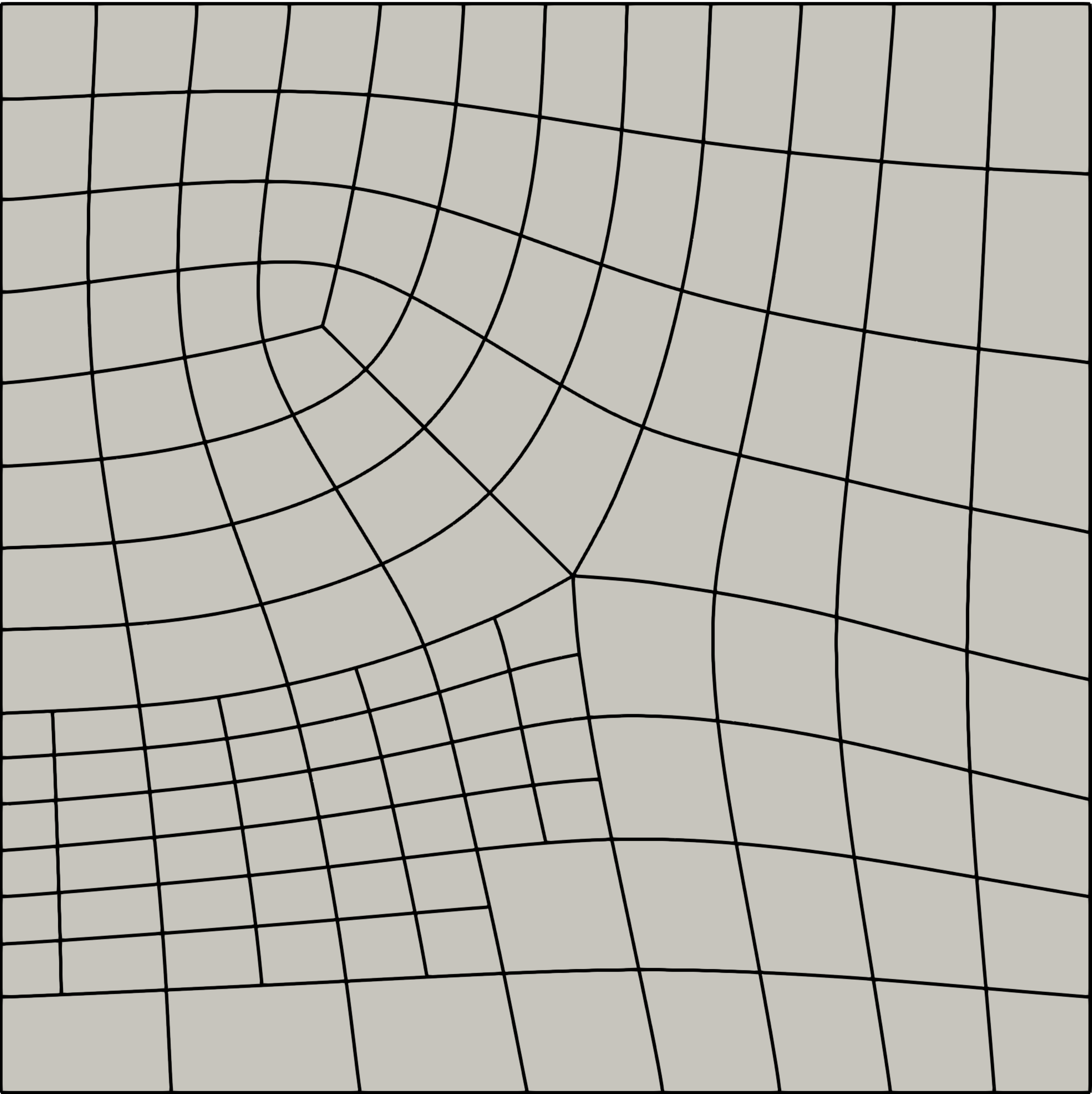} \\
(a) & (b) \\
\includegraphics[width=.47 \textwidth]{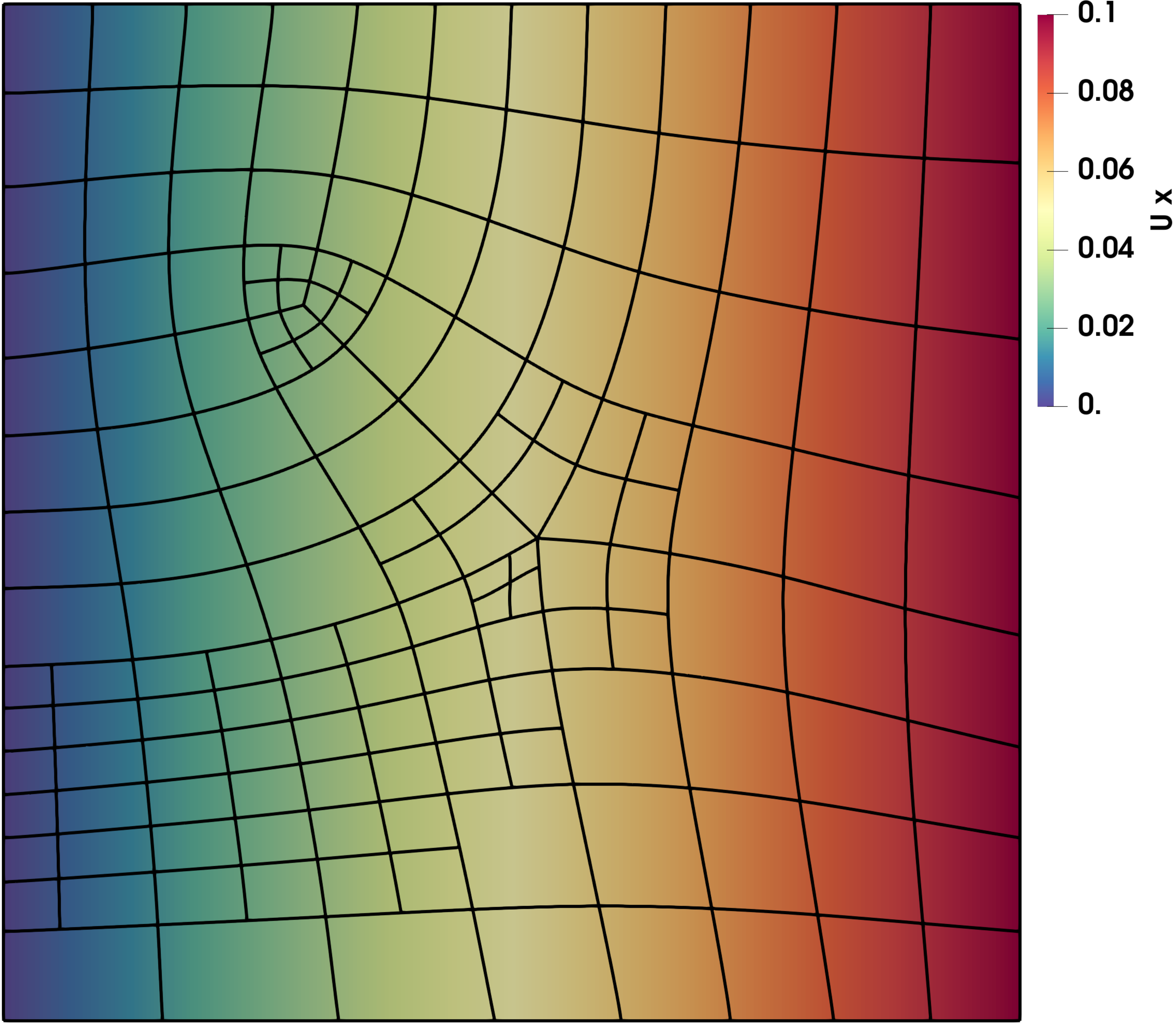} & \hspace{+1mm}
\includegraphics[width=.47 \textwidth]{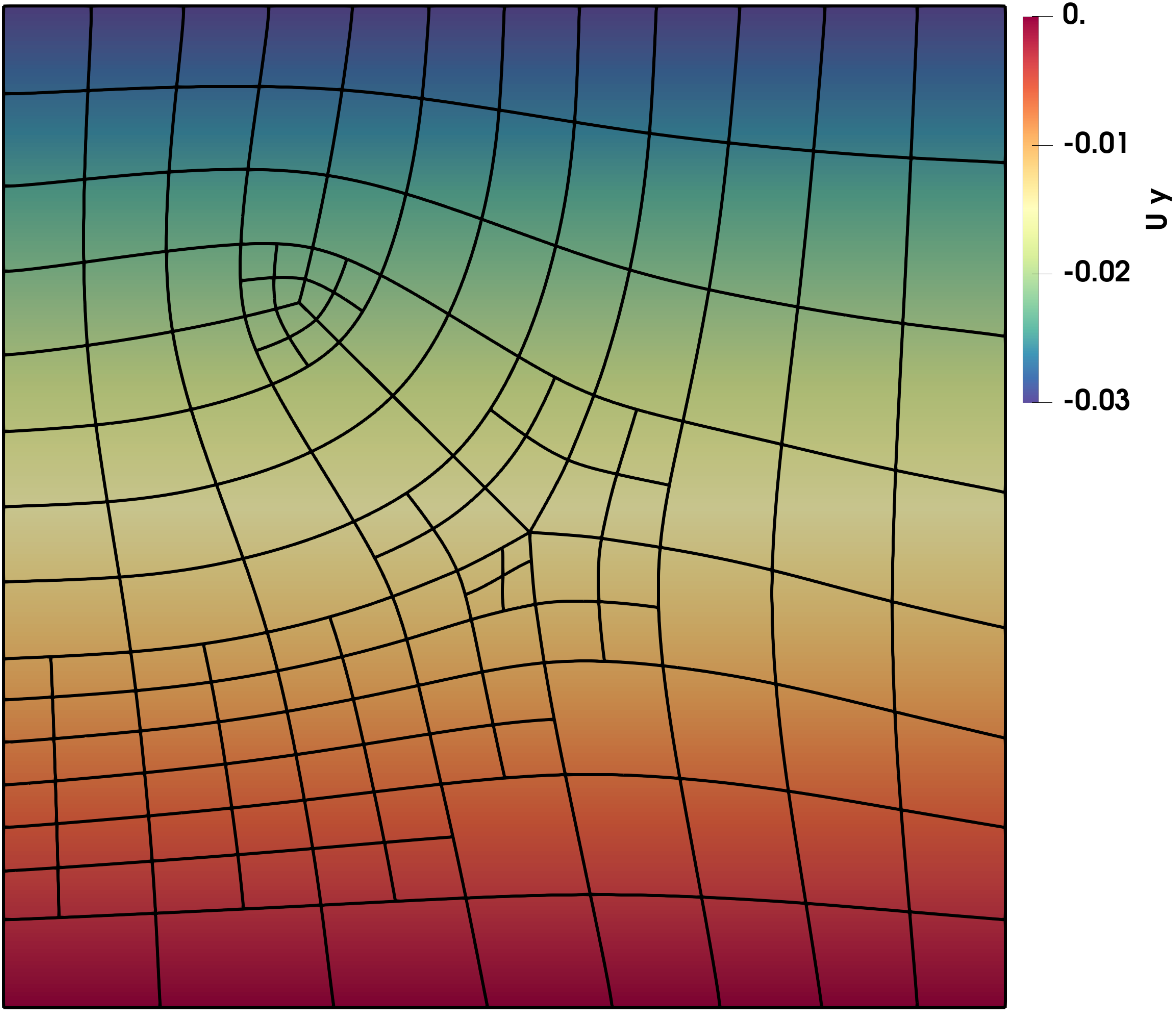} \\
(c) & (d) \\
\end{tabular}
\caption{The patch test using THU-splines on a unit square, where linear elasticity is solved. (a) The spline mesh of the initial local D-patch representation, whose input is given in Fig. \ref{fig:unstruct_mesh}; (b) the THU-spline mesh after several elements in (a) are refined; and (c, d) displacement ($u_x$ and $u_y$) results on the B\'{e}zier mesh. }
\label{fig:patch_test}
\end{figure}

In this section, we present several examples to show its effectiveness and efficiency. We start with a patch test on a domain of unit square, where linear elasticity is solved using THU-splines. The unstructured quad mesh in Fig. \ref{fig:unstruct_mesh} is taken as the input. Accordingly, an initial local D-patch representation is constructed, with its high-order spline mesh shown in Fig. \ref{fig:patch_test}(a). Next, we locally refine certain elements using THU-splines, which is intended to be general to include all types of elements: irregular, transition, and regular elements; see the resulting THU-spline mesh in Fig. \ref{fig:patch_test}(b). We use it to perform the patch test. Dirichlet boundary conditions $u_x=0$, $u_x=0.1$, and $u_y=0$ are imposed on the left, right, and lower boundaries, respectively, whereas free traction is applied to all the other boundaries. Isotropic and homogeneous material is assumed, where the Young's modulus $E$ and the Poisson's ratio are given as $E=1$ and $\nu=0.3$, respectively. The displacement results are shown on the B\'{e}zier mesh in Fig. \ref{fig:patch_test}(c, d), where recall that each irregular element in Fig. \ref{fig:patch_test}(b) needs to be split into $2\times 2$ B\'{e}zier elements for the smooth D-patch construction. Moreover, we find that the $L^2$-norm error of strains ($\epsilon_{xx}$, $\epsilon_{yy}$ and $\epsilon_{xy}$) is in the order of $10^{-16}$. In other words, the patch is passed with machine precision.

Successfully passing the patch test shows the applicability of THU-splines in engineering analysis. The numerical results also show supportive evidence regarding important properties of THU-splines. First, the involved THU-spline functions sum up to one and thus form a partition of unity. This confirms that the two kinds of truncation work as expected. Second, the involved THU-splines are linearly independent and thus constitute a basis. However, all these properties need to be rigorously studied, which is also the primary task of the future work.

\begin{figure}[htb]
\centering
\begin{tabular}{cc}
\includegraphics[width=.4 \textwidth]{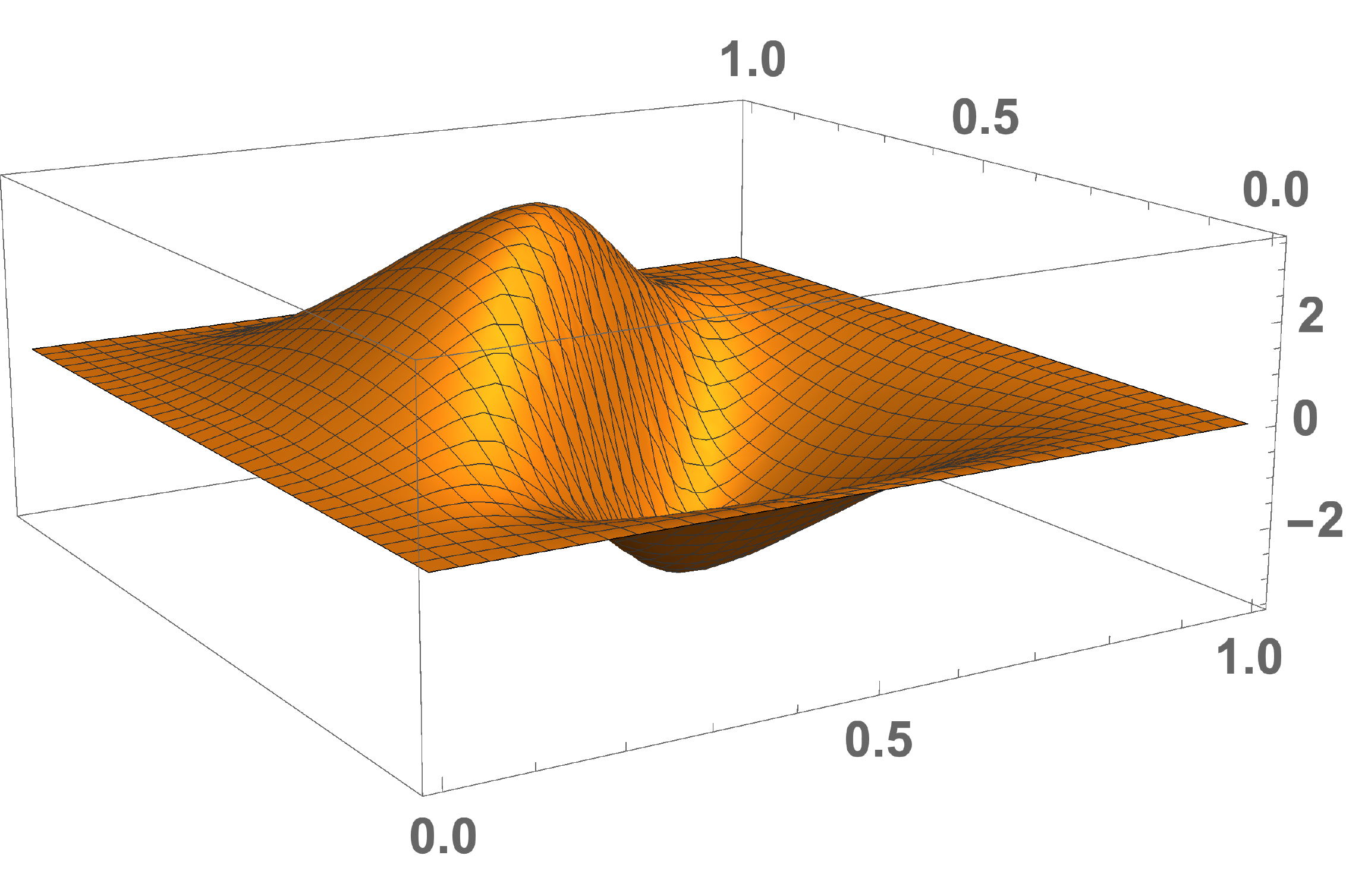} & \hspace{-4mm}
\includegraphics[width=.5 \textwidth]{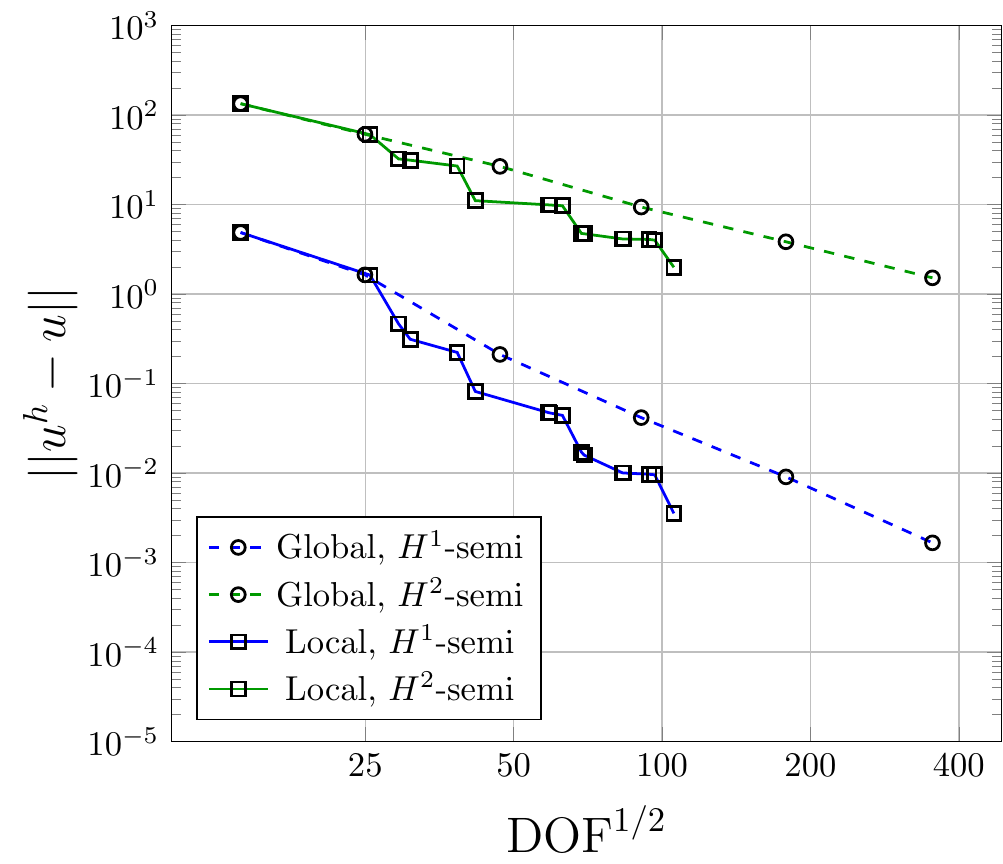} \\
(a) & (b)  \\
\includegraphics[height=.4 \textwidth]{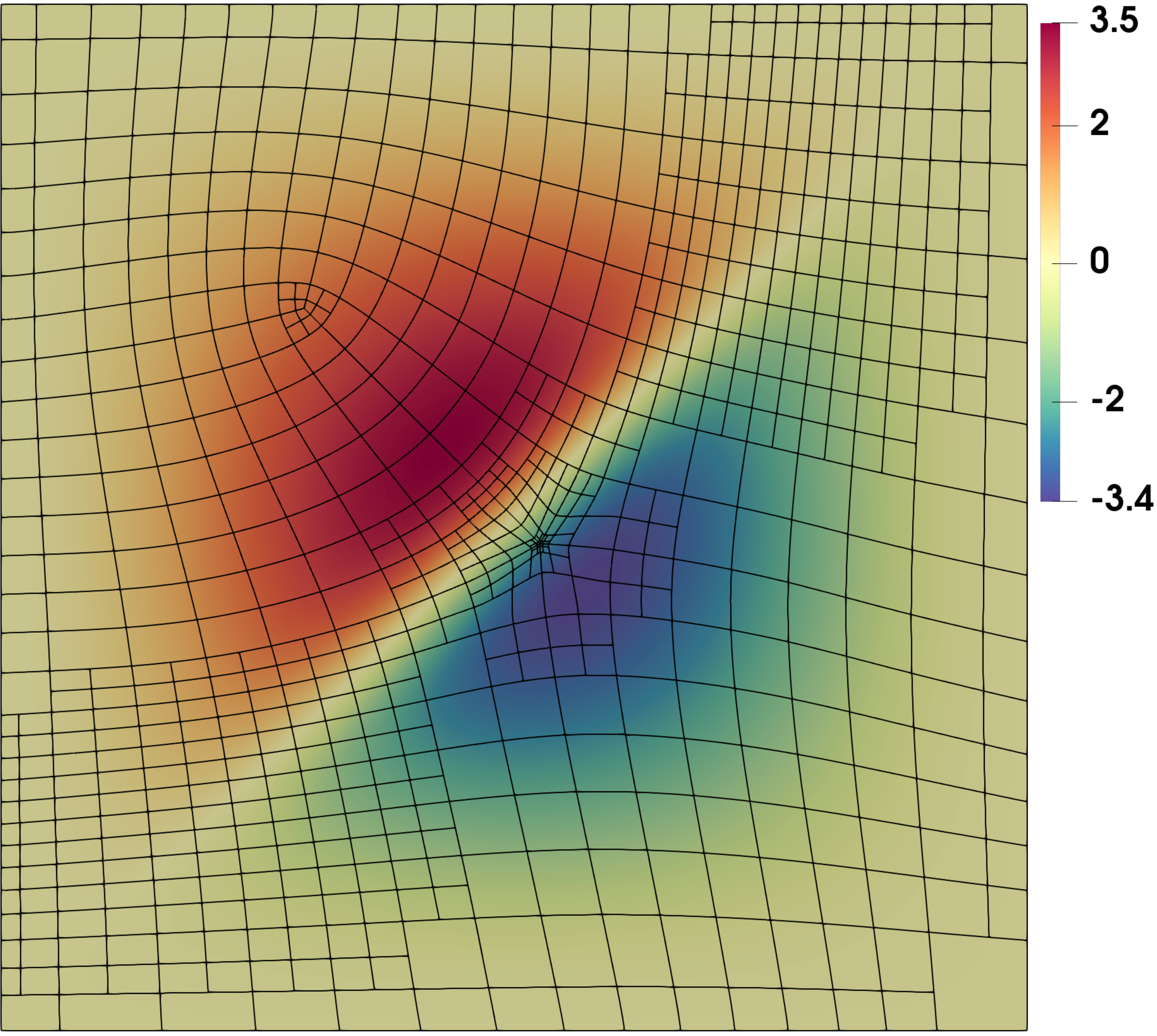} & \hspace{+2mm}
\includegraphics[height=.4 \textwidth]{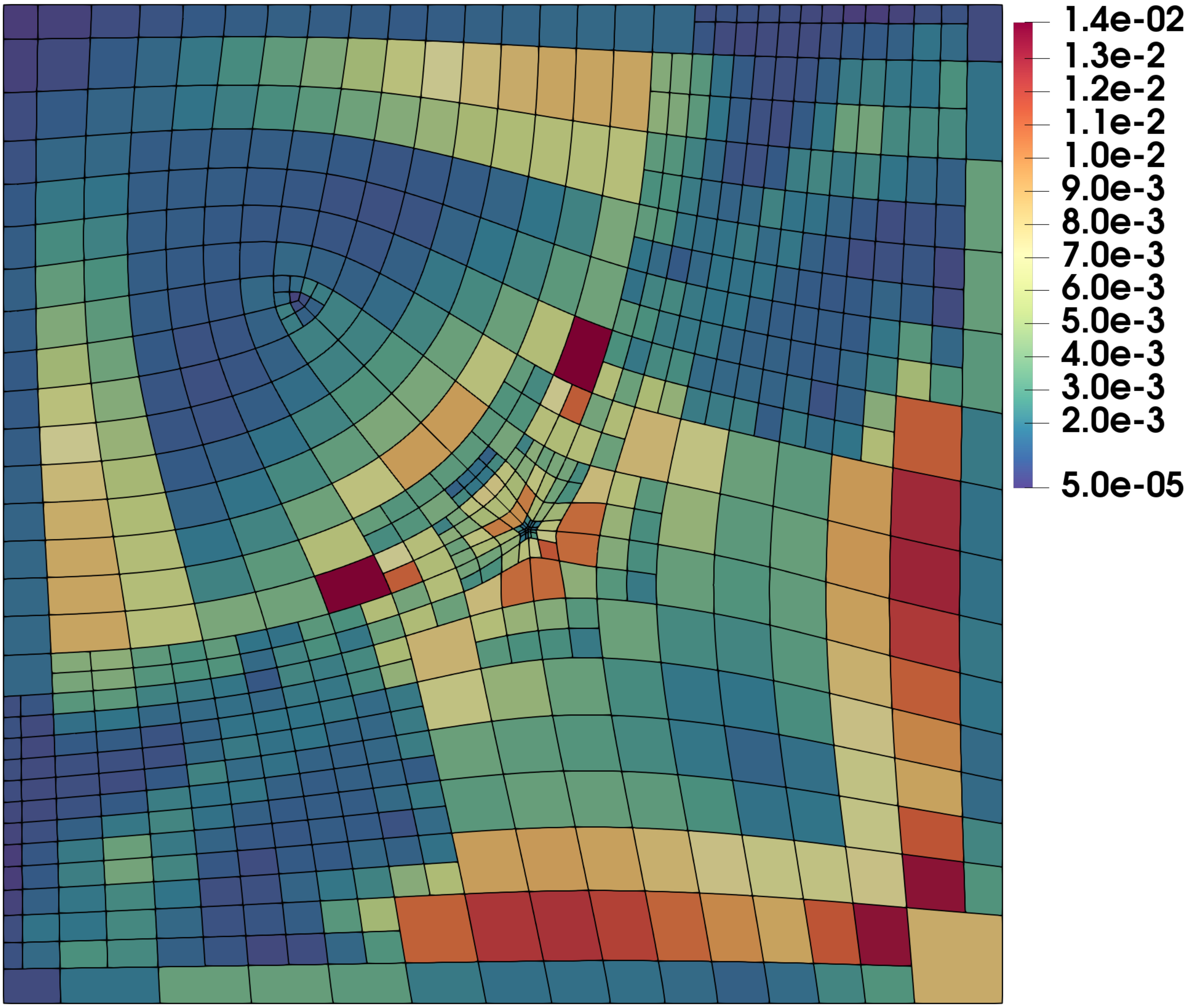} \\
(c)  & (d) \\
\end{tabular}
\caption{Adaptively solving the biharmonic equation using THU-splines. (a) The graph of the manufactured solution, (b) the convergence plot with respect to $\mathrm{DOF}^{\frac{1}{2}}$, (c) the approximate solution obtained on a locally refined mesh, and (d) the element-wise $H^1$ semi-norm error.}
\label{fig:biharmonic}
\end{figure}

Next, we use THU-splines to adaptively solve the biharmonic equation with the homogeneous Dirichlet and Neumann boundary conditions,
\begin{equation}
\left\{
\begin{array}{rll}
\Delta^2 u &= f \quad &\text{in} \quad \Omega, \\
 u &= 0 \quad &\text{on} \quad \partial\Omega, \\
 \nabla u \cdot \bm{n} &= 0 \quad &\text{on} \quad \partial\Omega, \\
\end{array}
\right.
\end{equation}
where $\bm{n}$ is the outwards normal of $\partial\Omega$.
This is a 4th order PDE that requires $C^1$ continuity of the underlying basis functions. We again solve the problem on a unit square domain and start with the mesh given in Fig. \ref{fig:patch_test}(a). For the sake of convergence study, we use the following manufactured solution that has a large gradient across the line $x-y=0$,
\begin{equation}
u(x,y) = a x^2 (1-x)^2 y^2 (1-y)^2 \tanh (b (y-x)), \quad (x,y) \in [0,1]^2,
\end{equation}
where $a$ and $b$ are two constants that control the behavior of the solution, and particularly, $b$ controls the ``sharpness" of the gradient; see Fig. \ref{fig:biharmonic}(a), where we use $a=10^3$ and $b=10\sqrt{2}$. $H^1$ semi-norm error is used to drive local refinement to efficiently capture the gradient feature. Global refinement is taken as a reference. The convergence plot is summarized in Fig. \ref{fig:biharmonic}(b) with respect to the square root of degrees of freedom (DOF). The approximate solution $u^h$ and the element-wise $H^1$ semi-norm error on the final locally refined mesh are shown in Fig. \ref{fig:biharmonic}(c) and (d), respectively. We observe that compared to global refinement, THU-splines can achieve the same accuracy with much fewer DOF.

\section{Conclusion and future work}
\label{sec:conclude}

In this paper, we have presented the construction method of THU-splines, a novel method that features a global smooth construction with the presence of EVs and supports flexible local refinement through the idea of THB-splines. Numerical evidence shows that THU-spline functions are refinable and form a partition of unity. Moreover, through solving the biharmonic equation with a manufactured solution, we show the efficiency of THU-splines compared to global refinement.

However, much remains to be studied to fully unleash the power of THU-splines. First of all, we will rigorously prove properties of THU-splines, including refinability, partition of unity, and linear independence. Once these foundations are ready, we will apply THU-splines to Kirchhoff-Love shell problems, where a posteriori error estimate is needed to automatically guide the adaptive analysis of THU-splines. The recent work reported in \cite{ref:coradello20} may well fit this purpose with the adaptation in irregular regions.

\begin{acknowledgement}
X. Wei is partially supported by the ERC AdG project CHANGE n. 694515, as well as the Swiss National Science Foundation project HOGAEMS n.200021\_188589.
\end{acknowledgement}

\bibliographystyle{plain}
\bibliography{ref}

\end{document}